\documentclass[12pt]{elsarticle}
\usepackage{amsfonts,enumerate,amsmath,amssymb,mathrsfs,amsbsy}
\usepackage[colorlinks=true]{hyperref}

\usepackage{caption}
\biboptions{numbers,sort&compress}

\usepackage{xcolor}

\def\qed{\strut\hfill $\Box$}
\newtheorem{thm}{Theorem}[section]

\newtheorem{lem}[thm]{Lemma}

\newtheorem{rem}[thm]{Remark}

\newcommand{\thmref}[1]{Theorem~{\rm \ref{#1}}}
\newcommand{\lemref}[1]{Lemma~{\rm \ref{#1}}}

\newcommand{\remref}[1]{Remark~{\rm \ref{#1}}}

\def\para#1{\vskip .4\baselineskip\noindent{\bf #1}}
\numberwithin{equation}{section}
\allowdisplaybreaks \allowdisplaybreaks[4]
\begin{document}
\begin{frontmatter}	
	\title{Averaging principles for non-autonomous two-time-scale stochastic
		reaction-diffusion equations with jump}
	
	\author[mymainaddress,myfivearyaddress]{Yong Xu\corref{mycorrespondingauthor}}
	
	\cortext[mycorrespondingauthor]{Corresponding author}
	\ead{hsux3@nwpu.edu.cn}
	
	\author[mymainaddress]{Ruifang Wang}
	\ead{wrfjy@yahoo.com}

	
	
	
	\address[mymainaddress]{Department of Applied Mathematics, Northwestern Polytechnical University, Xi'an, 710072, China}
	\address[myfivearyaddress]{MIIT Key Laboratory of Dynamics and Control of Complex Systems, Northwestern Polytechnical University, Xi'an, 710072, China	}
	

	%

	\begin{abstract}
		In this paper, we aim to develop the averaging principle for a slow-fast system of stochastic reaction-diffusion equations driven by Poisson random measures. The coefficients of the equation are assumed to be functions of time, and some of them are periodic or almost periodic. Therefore, the Poisson term needs to be processed, and a new averaged equation needs to be given. For this reason, the existence of time-dependent evolution family of measures associated with the fast equation is studied, and proved that it is almost periodic. Next, according to the characteristics of  almost periodic functions, the averaged coefficient is defined by the evolution family of measures, and the averaged equation is given. Finally, the validity of the averaging  principle is verified by using the Khasminskii method.	
	\end{abstract}		
\end{frontmatter}

\section{Introduction}\label{sec-1}
The slow-fast systems are widely encountered in biology, ecology and other application areas.   
In this paper, we are concerned with the following non-autonomous slow-fast systems of stochastic partial differential equations (SPDEs) on a bounded domain $\mathcal{O} $  of  $ \mathbb{R}^d\left( d\ge 1 \right)  $: 
\begin{eqnarray}\label{orginal1}
\left\{\begin{array}{l} 
\frac{\partial u_{\epsilon}}{\partial t}\left( t,\xi \right) =\mathcal{A}_1\left(t \right) u_{\epsilon}\left( t,\xi \right) +b_1\left(t, \xi ,u_{\epsilon}\left( t,\xi \right),v_{\epsilon}\left( t,\xi \right) \right)\\
\qquad\qquad\quad+f_1\left(t, \xi ,u_{\epsilon}\left( t,\xi \right) \right) \frac{\partial \omega ^{Q_1}}{\partial t}\left( t,\xi \right) \\
\qquad\qquad\quad+\int_{\mathbb{Z}}{g_1\left(t, \xi ,u_{\epsilon}\left( t,\xi \right) ,z \right) \frac{\partial \tilde{N}_1}{\partial t}\left( t,\xi ,dz \right)},\\
\frac{\partial v_{\epsilon}}{\partial t}\left( t,\xi \right) 
= \frac{1}{\epsilon}\left[ \left( \mathcal{A}_2\left( t \right) -\alpha \right) v_{\epsilon}\left( t,\xi \right) +b_2\left( t,\xi ,u_{\epsilon}\left( t,\xi \right) ,v_{\epsilon}\left( t,\xi \right) \right) \right] \\
\qquad\qquad\quad+\frac{1}{\sqrt{\epsilon}}f_2\left( t,\xi ,u_{\epsilon}\left( t,\xi \right) ,v_{\epsilon}\left( t,\xi \right) \right) \frac{\partial \omega ^{Q_2}}{\partial t}\left( t,\xi \right) \\
\qquad\qquad\quad+\int_{\mathbb{Z}}{g_2}\left( t,\xi ,u_{\epsilon}\left( t,\xi \right) ,v_{\epsilon}\left( t,\xi \right) ,z \right) \frac{\partial \tilde{N}_{2}^{\epsilon}}{\partial t}\left( t,\xi ,dz \right) ,\\
u_{\epsilon}\left( 0,\xi \right) =x\left( \xi \right), \quad 
v_{\epsilon}\left( 0,\xi \right) =y\left( \xi \right), \quad
\xi \in\mathcal{O},\\
\mathcal{N}_1u_{\epsilon}\left( t,\xi \right) = \mathcal{N}_2v_{\epsilon}\left( t,\xi \right) =0, \quad t\ge 0, \quad\xi \in \partial\mathcal{O},
\end{array}
\right.
\end{eqnarray}
where $ \epsilon \ll 1 $ is a  positive parameter and $ \alpha $  is a sufficiently large fixed constant. The operators $\mathcal{N}_1$ and $\mathcal{N}_2$ are boundary operators. The stochastic perturbations $ \omega ^{Q_1}, \omega ^{Q_2} $ and $ \tilde{N}_1, \tilde{N}_{2}^{\epsilon} $ are mutually independent Wiener processes and Poisson random measures on the same complete stochastic basis $\big( \varOmega ,\mathcal{F},\left\lbrace \mathcal{F}_t\right\rbrace _{t\geq0},\mathbb{P} \big), $ which  will be described in Section \ref{sec-2}. For $ i=1,2 $, the operator $ \mathcal{A}_i\left( t\right)  $  and the functions $ b_i, f_i, g_i$  depend on time, and we assume that the operator $ \mathcal{A}_2\left(t \right)  $ is periodic and the functions $ b_1, b_2, f_2, g_2$ are almost periodic.

The goal of this paper is to establish an effective approximations for the slow equation of the original system (\ref{orginal1})  by using the averaging principle. The averaged equation is obtained as following
\begin{eqnarray}\label{en15}
\begin{split}
\begin{cases}
\frac{\partial \bar{u}}{\partial t}\left( t,\xi \right) &=\mathcal{A}_1\left(t \right) \bar{u}\left( t,\xi \right) +\bar{B}_1\left( \bar{u}\left( t \right) \right)\left( \xi\right)  +f_1\left(t, \xi ,\bar{u}\left( t,\xi \right) \right) \frac{\partial \omega ^{Q_1}}{\partial t}\left( t,\xi \right) \\
&\quad+\int_{\mathbb{Z}}{g_1\left(t, \xi ,\bar{u}\left( t,\xi \right) ,z \right) \frac{\partial \tilde{N}_1}{\partial t}\left( t,\xi ,dz \right)},\\
\bar{u}\left( 0,\xi \right) &=x\left( \xi \right), \quad \xi\in\mathcal{O}, \quad
\mathcal{N}_1\bar{u}\left( t,\xi \right) =0, \quad t\ge 0,\quad \xi \in \partial\mathcal{O},
\end{cases} 
\end{split}
\end{eqnarray}
where $ \bar{B}_1 $ is the averaged coefficient, which will be given in equation (\ref{en18}). To demonstrate the validity of the averaging principle, we prove that for any  $ T>0 $ and  $ \eta >0 $, it yields 
\begin{eqnarray}\label{en711}
\underset{\epsilon \rightarrow 0}{\lim}\mathbb{P}\Big( \underset{t\in \left[ 0,T \right]}{\sup}\left\| u_{\epsilon}\left( t \right) -\bar{u}\left( t \right) \right\| _{L^2(\mathcal{O})}>\eta \Big) =0,
\end{eqnarray}
where  $ \bar{u} $ is the solution of the averaged equation (\ref{en15}).

The theory of  the averaging principle 
originated by Laplace and Lagrange, has been applied in celestial mechanics, oscillation theory, radiophysics and other fields. 
The firstly rigorous results for the deterministic case were given by Bogolyubov and  Mitropolskii \cite{Bogolyubov1961asymptotic}.  Moreover, Volosov \cite{Volosov1962Averaging} and Besjes \cite{besjes1969asymptotic} also promoted the development of the averaging principle. Then, great interests have appeared in its application to dynamical systems under random perturbations. An important contribution was that, in 1968,  Khasminskii \cite{khas1968on} originally proposed the averaging principle for stochastic differential equations (SDEs) driven by Brownian motion. Since then, the averaging principle has been an active area of research. Many studies on the averaging principle of SDEs have been presented, e.g, Givon \cite{givon2007strong}, Freidlin and Wentzell \cite{freidlin2012random},  Duan \cite{duan2014effective}, Thompson \cite{thompson2015stochastic}, Xu and his co-workers
\cite{xu2011averaging,xu2015stochastic,xu2017stochastic}. Recently,  effective approximations for slow-fast  SPDEs have been received extensive attention. Cerrai \cite{cerrai2009khasminskii,cerrai2011averaging} investigated the validity of the averaging principle for a class of stochastic reaction-diffusion equations with multiplicative noise. In addition, Wang and Roberts \cite{wang2012average}, Pei and Xu   \cite{pei2017averaging,pei2020Stochastic,pei2020convergence}, Li \cite{li2019stochastic} and Gao \cite{gao2018averaging,gao2018kuramot,gao2019averaging} also studied the averaging principles for the slow-fast SPDEs.

The above-mentioned works mainly considered autonomous systems. For  the autonomous systems, as long as the initial value is given, the solution of which only depends on the duration of time, not on the selection of the initial time. However, if the initial time is different, the solution of non-autonomous equations with the same initial data will also be different. Therefore, compared with  the autonomous systems, the dynamic behaviors of the non-autonomous systems  are more complex, which can portray more actual models. 
Chepyzhov and Vishik \cite{chepyzhov2002attractors} studied  long-time dynamic behaviors of non-autonomous dissipative system.
Carvalho \cite{Carvalho2012attractors} dealt with the theory of attractors for  the non-autonomous dynamical systems. Bunder and Roberts \cite{Bunder2017resolution} considered the discrete modelling of non-autonomous PDEs. 

In 2017, the averaging principle has been presented for non-autonomous slow-fast system of stochastic reaction-diffusion equations by Cerrai \cite{cerrai2017averaging}. But, the system of this paper was driven by Gaussian noises, which is considered as an ideal noise source and only can simulate fluctuations near the mean value. Actually, due to the complexity of the external environment, random noise sources encountered in practical fields usually exhibit non-Gaussian properties, which may cause sharply fluctuations. It should be pointed out that Poisson noise, one of the most ubiquitous noise sources in many fields \cite{Applebaum2009Processes,duan2015introduction,peszat2007stochastic}, can provide a accurate mathematical model to describe discontinuous random processes, some large moves and unpredictable events \cite{Bertoin1998levy,zhang2019random,li2018stepanov-like}. So, in this paper, we are devoted to developing the averaging principle for non-autonomous systems of reaction-diffusion equations driven by  Wiener processes and Poisson random measures. 

The key to using the averaging principle to analyze system  (\ref{orginal1}) is the fast equation with a frozen slow component $x\in L^2\left(\mathcal{O}\right) $:
\begin{eqnarray}\label{en11}
\begin{split}
\begin{cases}
\frac{\partial v^{x,y}}{\partial t}\left( t,\xi \right) 
&= \left[ \left( \mathcal{A}_2\left( t \right) -\alpha \right) v^{x,y}\left( t,\xi \right) +b_2\left( t,\xi ,x\left( \xi \right) ,v^{x,y}\left( t,\xi \right) \right) \right] \\
&\quad+f_2\left( t,\xi ,x\left( \xi \right),v^{x,y}\left( t,\xi \right) \right) \frac{\partial \omega ^{Q_2}}{\partial t}\left( t,\xi \right) \\
&\quad+\int_{\mathbb{Z}}{g_2}\left( t,\xi,x\left( \xi \right),v^{x,y}\left( t,\xi \right) ,z \right) \frac{\partial \tilde{N}_{2}^{\epsilon} }{\partial t}\left( t,\xi ,dz \right) ,\\
v^{x,y}\left( s,\xi \right) &=y\left( \xi \right), \quad
\xi \in\mathcal{O}, \qquad \mathcal{N}_2v^{x,y}\left( t,\xi \right) =0, \quad t\ge s, \quad\xi \in \partial\mathcal{O}.
\end{cases} 
\end{split}
\end{eqnarray}


By dealing with the Poisson terms, we prove that an evolution family of measures  $ \left( \mu _{t}^{x};t\in \mathbb{R} \right) $ on $ L^2\left( \mathcal{O}\right) $ for the fast equation (\ref{en11}) exists. Then, assuming that $  \mathcal{A}_2(t) $ is periodic and $ b_2, f_2, g_2 $ are almost periodic, we prove that the evolution family of measures is almost periodic.
Moreover, with the aid of the Theorem 2.10 in \cite{fink1974almost}, we prove that the family of functions 
\begin{eqnarray}\label{en14}
\left\{ t\in \mathbb{R}\mapsto \int_{L^2\left(\mathcal{O}\right)}{B_1\left(t, x,y \right)}\mu _{t}^{x}\left( dy \right) \right\}  \nonumber
\end{eqnarray}
is uniformly almost periodic for any $ x $ in compact set $ \mathbb{K}\subset L^2\left(\mathcal{O}\right) $, where $B_1(t, x, y)(\xi )=b_1(t, \xi, x( \xi ), y( \xi ) )$ for any $ x, y \in L^2\left(\mathcal{O}\right) $ and $ \xi\in\mathcal{O} $. 

According to the characteristics of almost periodic function \cite[Theorem 3.4]{cerrai2017averaging}, we  define the averaged coefficient $ \bar{B}_1 $  
as follows
\begin{eqnarray}\label{en18}
\bar{B}_1\left( x \right) :=\underset{T\rightarrow \infty}{\lim}\frac{1}{T}\int_0^T{\int_{L^2\left(\mathcal{O}\right)}{B_1\left( t,x,y \right)}\mu _{t}^{x}\left( dy \right) dt, \quad x\in L^2\left(\mathcal{O}\right)}.
\end{eqnarray}

Finally, the averaged equation is obtained through the averaged coefficient $ \bar{B}_1 $. Using the classical Khasminskii method to the present situation, the averaging principle is effective.

The above-mentioned notations will be given in Section \ref{sec-2}, and in this paper, $ c>0 $  below with or without subscripts will represent a universal constant whose value may vary in different occasions.

\section{Notations, assumptions and preliminaries}\label{sec-2}
Let $\mathcal{O}$ be a bounded domain of $\mathbb{R}^d\left( d\ge 1 \right)$ having a  smooth boundary. In this paper, we denote $\mathbb{H}$  the separable Hilbert space $L^2\left(\mathcal{O} \right)$, endowed with the usual scalar product 
$$\left< x , y \right> _{\mathbb{H}}=\int_{\mathcal{O}} x\left(\xi \right)y\left(\xi \right) d\xi $$ 
and with the corresponding norm $\left\| \cdot \right\| _{\mathbb{H}}.$ 
The norm in $ L^\infty \left( \mathcal{O}\right)  $ will be denoted by $ \left\| \cdot\right\|_\infty $.

Furthermore, the subspace $ \mathcal{D}(( -A ) ^{\theta})  $ \cite{boufoussi2010successive,luo2008fixed,luo2009existence} of the generator $ A $ is dense in $ \mathbb{H} $, and endowed with the norm
$$
\left\| \varLambda \right\| _{\theta}=\lVert \left( -A \right) ^{\theta}\varLambda \rVert_{\mathbb{H}} , \quad \varLambda \in \mathcal{D}(( -A ) ^{\theta}),
$$
for $ 0 \leq \theta < 1,0<t\le T $. According to Theorem 6.13 in \cite{pazy1983semigroups}, there exists a $ c_{\theta}>0 $, such that
$$
\lVert \left( -A \right) ^{\theta}e^{At} \rVert_{\mathbb{H}} \le c_{\theta}t^{-\theta}.
$$

Denote by $B_b\left( \mathbb{H} \right)$ the Banach space of the bounded Borel functions $\varphi :\mathbb{H}\rightarrow \mathbb{R}$, endowed with the sup-norm 
$$
\left\| \varphi \right\| _0:=\underset{x\in \mathbb{H}}{\sup}| \varphi \left( x \right)|,
$$
and $ C_b(\mathbb{H}) $ is the subspace of the uniformly continuous mappings. Moreover, $  D\left( \left[ s,T \right] ;  \mathbb{H}  \right) $ denotes the space of all c\`adl\`ag path from $ \left[ s,T \right] $ into $  \mathbb{H}.  $ 

We shall denote that $\mathcal{L}\left( \mathbb{H} \right)$ is the space of the bounded linear operators in $\mathbb{H}$, and denote $\mathcal{L}_2\left( \mathbb{H} \right)$ the subspace of Hilbert-Schmidt operators, endowed with the norm $$
\left\| Q \right\| _2=\sqrt{Tr\left[ Q^{\ast}Q \right]}.
$$

In the slow-fast system (\ref{orginal1}), the Gaussian noises ${\partial \omega ^{Q_1}}/{\partial t}\left( t,\xi \right)$ and ${\partial \omega ^{Q_2}}/{\partial t}\left( t,\xi \right)$ are assumed to be white in time and colored in space in the case of space dimension $d>1$, for $t\ge 0$ and $\xi \in\mathcal{O}$. And, $\omega ^{Q_i}\left( t,\xi \right) \left( i=1,2 \right)$ is the cylindrical Wiener processes, defined as $$
\omega ^{Q_i}\left( t,\xi \right) =\sum_{k=1}^{\infty}{Q_ie_k\left( \xi \right) \beta _k\left( t \right)}, \quad i=1,2,
$$
where $\left\{ e_k \right\} _{k\in \mathbb{N}}$ is a complete orthonormal basis in $\mathbb{H}$, $\left\{ \beta _k\left( t \right) \right\} _{k\in \mathbb{N}}$ is a sequence of mutually independent standard Brownian motion defined on the same complete stochastic basis $\big( \varOmega ,\mathcal{F},\left\lbrace \mathcal{F}_t\right\rbrace _{t\geq0},\mathbb{P} \big) $,  $Q_i$ is a bounded linear operator on $\mathbb{H}$.

Next, we give the definitions of Poisson random measures $\tilde{N}_1\left( dt,dz \right)$ and $\tilde{N}_{2}^{\epsilon}\left( dt,dz \right)$. Let $ \left( \mathbb{Z},\mathcal{B}\left( \mathbb{Z} \right) \right) $  be a given measurable space and  $ v\left( dz \right)  $ be a  $ \sigma  $-finite measure on it. $ D_{p_{t}^{i}},i=1,2 $ are two countable subsets of $ \mathbb{R}_+ $. Moreover, let
$ p_{t}^{1},t \in  D_{p_{t}^{1}}  $ be a stationary $ \mathcal{F}_t $-adapted Poisson point process on $ \mathbb{Z} $ with the characteristic $ v $,
and $ p_{t}^{2},t \in  D_{p_{t}^{2}} $ be the other stationary $ \mathcal{F}_t $-adapted Poisson point process on $ \mathbb{Z} $ with the characteristic $ {v}/{\epsilon} $. Denote by $ N_i\left( dt,dz \right),i=1,2 $ the Poisson counting measure associated with $ p_{t}^{i} $, i.e.,
$$
N_i\left( t,\varLambda \right) :=\sum_{s\in D_{p_{t}^{i}},s\le t}{I_{\varLambda}\left( p_{t}^{i} \right)},\quad i=1,2.
$$
Let us denote the two independent compensated Poisson measures
$$\tilde{N}_1\left( dt,dz \right) :=N_1\left( dt,dz \right) -v_1\left( dz \right) dt$$ 
and 
$$\tilde{N}_{2}^{\epsilon}\left( dt,dz \right) :=N_2\left( dt,dz \right) -\frac{1}{\epsilon}v_2\left( dz \right) dt,$$ 
where  $v_1\left( dz \right) dt$ and $\frac{1}{\epsilon}v_2\left( dz \right) dt$ are the compensators. 

Refer to \cite{prato2014stochastic,peszat2007stochastic} for a more detailed description of the stochastic integral with respect to a cylindrical Wiener process and Poisson random measure. 

For any $ t\in \mathbb{R} $, the operators $\mathcal{A}_1\left( t \right) $ and $\mathcal{A}_2\left( t \right) $ are second order uniformly elliptic operators, having continuous coefficients on $\mathcal{O}$. The operators $\mathcal{N}_1$ and $\mathcal{N}_2$ are the boundary operators, which can be either the identity operator (Dirichlet boundary condition) or a first order operator (coefficients satisfying a uniform nontangentiality condition). We shall assume that the operator $ \mathcal{A}_i(t) $ has the following form
\begin{eqnarray}\label{en21}
\begin{split}
\mathcal{A}_i\left( t \right) =\gamma_i \left( t \right) \mathcal{A}_i+\mathcal{L}_i\left( t \right), \quad t\in \mathbb{R},\ i=1,2,
\end{split}
\end{eqnarray}
where $\mathcal{A}_i$ is a second order uniformly elliptic operator \cite{Evans2010Partial,Amann1995Linear} with continuous coefficients on $\mathcal{O}$, which is independent of $t$. In addition, $\mathcal{L}_i\left( t \right) $ is a first order differential operator, has the form
\begin{eqnarray}\label{en22}
\begin{split}
\mathcal{L}_i\left( t,\xi \right) u\left( \xi \right) =\left< l_i\left( t,\xi \right) ,\nabla u\left( \xi \right) \right> _{\mathbb{R}^d}, \quad t\in \mathbb{R}, \ \xi \in \mathcal{O}.
\end{split}
\end{eqnarray}
The realizations of the differential operators $ \mathcal{A}_i $ and $ \mathcal{L}_i $ in $ \mathbb{H} $ is $ A_i $ and $ L_i $. Moveover, $ {A}_1 $ and $ {A}_2 $ generate two analytic semigroups $ e^{tA_1} $ and $ e^{tA_2} $ respectively. 

Now, we give the following assumptions:	
\begin{enumerate}[({A}1)]
	\item
	\begin{enumerate}
		\item 	For $ i=1,2 $, the function $\gamma_i :\mathbb{R}\rightarrow \mathbb{R}$ is continuous, and there exist $\gamma _0, \gamma >0$ such that 
		\begin{eqnarray}\label{en23}
		\begin{split}
		\gamma _0\le \gamma_i \left( t \right) \le \gamma, \quad t\in \mathbb{R}.
		\end{split}
		\end{eqnarray}
		\item For $ i=1,2 $, the function $l_i:\mathbb{R}\times \mathcal{O}\rightarrow \mathbb{R}^d$ is continuous and bounded.
	\end{enumerate}	
	\item For $i=1,2$, there exist a complete orthonormal system $\left\{ e_{i,k} \right\} _{k\in \mathbb{N}}$ in $\mathbb{H}$ and two sequences of nonnegative real numbers $\left\{ \alpha _{i,k} \right\} _{k\in \mathbb{N}}$ and $\left\{ \lambda _{i,k} \right\} _{k\in \mathbb{N}}$ such that
	\begin{eqnarray}\label{en24}
	\begin{split}
	A_ie_{i,k}=-\alpha _{i,k}e_{i,k}, \quad Q_ie_{i,k}=\lambda _{i,k}e_{i,k}, \quad k\ge 1,
	\end{split}
	\end{eqnarray}
	and
	\begin{eqnarray}\label{en25}
	\begin{split}
	\kappa _i:=\sum_{k=1}^{\infty}{\lambda _{i,k}^{\rho _i}\left\| e_{i,k} \right\| _{\infty}^{2}}<\infty, \quad \zeta _i:=\sum_{k=1}^{\infty}{\alpha _{i,k}^{-\beta _i}\left\| e_{i,k} \right\| _{\infty}^{2}}<\infty,
	\end{split}
	\end{eqnarray}
	for some constants $\rho _i\in \left( 2,+\infty \right]$ and $\beta _i\in \left( 0,+\infty \right)$ such that
	\begin{eqnarray}\label{en26}
	[{\beta _i\left( \rho _i-2 \right)}]/{\rho _i}<1.
	\end{eqnarray}
\end{enumerate} 
\begin{rem} 
	{\rm For more comments and examples about the assumption (A2) of the operators $ A_i $ and $ Q_i $, reader can read \cite{cerrai2009khasminskii}.}
\end{rem} 
\begin{enumerate}[({A}3)]
	\item The mappings $b_1:\mathbb{R}\times\mathcal{O}\times \mathbb{R}^2\rightarrow \mathbb{R}, f_1:\mathbb{R}\times\mathcal{O}\times \mathbb{R}\rightarrow \mathbb{R}, g_1:\mathbb{R}\times\mathcal{O}\times \mathbb{R}\times \mathbb{Z}\rightarrow \mathbb{R}$ are measurable, and the mappings $b_1\left(t, \xi ,\cdot \right):\mathbb{R}^2\rightarrow \mathbb{R}, f_1\left( t,\xi ,\cdot \right):\mathbb{R}\rightarrow \mathbb{R}$,   $g_1\left( t,\xi ,\cdot ,z \right):\mathbb{R}\rightarrow \mathbb{R}$ are Lipschitz continuous and linearly growing, uniformly with respect to $(t,\xi,z) \in \mathbb{R}\times\mathcal{O}\times \mathbb{Z}$. Moreover, for all  $ p\ge 1 $, there exist positive constants $ c_1,c_2 $,  such that for all $  x_1,x_2  \in \mathbb{R}  $, we have
	$$
	\underset{\left(t, \xi\right)  \in   \mathbb{R}\times\mathcal{O} }{\text{sup}}\int_{\mathbb{Z}}{ \left|  g_1\left(t, \xi ,x_1,z \right) \right|^p } \upsilon _1\left( dz \right) \leq c_1\left( 1+\left| x_1\right|^p  \right) , 
	$$
	$$
	\underset{\left(t, \xi\right)  \in   \mathbb{R}\times\mathcal{O} }{\text{sup}}\int_{\mathbb{Z}}{\left|  g_1\left(t, \xi ,x_1,z \right) -g_1\left(t, \xi ,x_2,z \right)\right|  }^{p}\upsilon _1\left( dz \right) \leq c_2\left|  x_1-x_2\right| ^{p}. 
	$$
	\end{enumerate} 
	\begin{enumerate}[({A}4)]
	\item 
	The mappings $b_2:\mathbb{R}\times\mathcal{O}\times \mathbb{R}^2\rightarrow \mathbb{R}, f_2:\mathbb{R}\times\mathcal{O}\times \mathbb{R}^2\rightarrow \mathbb{R}, g_2:\mathbb{R}\times\mathcal{O}\times \mathbb{R}^2\times \mathbb{Z}\rightarrow \mathbb{R}$ are measurable, and the mappings $b_2\left( t,\xi ,\cdot \right):\mathbb{R}^2\rightarrow \mathbb{R},  f_2\left(t, \xi ,\cdot \right):\mathbb{R}^2\rightarrow \mathbb{R}, g_2\left( t,\xi,\cdot,z \right):\mathbb{R}^2\rightarrow \mathbb{R}$ are Lipschitz continuous and linearly growing, uniformly with respect to $(t,\xi,z) \in \mathbb{R}\times\mathcal{O} \times \mathbb{Z}$. Moreover, for all  $ q\ge 1 $, there exist positive constants $c_3,c_4 $,  such that for all $ (x_i,y_i) \in \mathbb{R}^2,i=1,2 $, we have
	$$
	\underset{\left(t, \xi\right)  \in   \mathbb{R}\times\mathcal{O} }{\text{sup}}\int_{\mathbb{Z}}{\left|  g_2\left( t,\xi ,x_1,y_1,z \right)\right| } ^{q}\upsilon _2\left( dz \right) \le c_3\left( 1+\left| x_1\right|  ^{q}+\left|  y_1\right|^{q} \right), 
	$$
	\begin{eqnarray}
	\underset{\left(t, \xi\right)  \in   \mathbb{R}\times\mathcal{O} }{\text{sup}}\int_{\mathbb{Z}}{\left|  g_2\left( t,\xi ,x_1,y_1,z \right) -g_2\left( t,\xi ,x_2,y_2,z \right)\right| }^{q}\upsilon _2\left( dz \right) &\le& c_4\left( \left|  x_1-x_2\right| ^{q}\right. \cr
	&&\quad \left. +\left| y_1-y_2\right|^{q} \right). \nonumber
	\end{eqnarray}
\end{enumerate}	
\begin{rem}\label{rem2.1}
	{\rm For any $\left( t,\xi \right) \in \mathbb{R}\times\mathcal{O}$ and $x,y,h\in \mathbb{H},z\in \mathbb{Z}$, we shall set
		\begin{gather}
		B_1\left(t, x,y \right) \left( \xi \right) :=b_1\left(t, \xi ,x\left( \xi \right) ,y\left( \xi \right) \right) , \quad B_2\left( t,x,y \right) \left( \xi \right) :=b_2\left( t,\xi ,x\left( \xi \right) ,y\left( \xi \right) \right) ,\cr
		\left[ F_1\left(t, x \right) h \right] \left( \xi \right) :=f_1\left(t, \xi ,x\left( \xi \right) \right) h\left( \xi \right), \quad \left[ F_2\left( t,x,y \right) h \right] \left( \xi \right) :=f_2\left( t,\xi ,x\left( \xi \right) ,y\left( \xi \right) \right) h\left( \xi \right) ,\cr
		\left[ G_1\left(t, x,z \right) h \right] \left( \xi \right) :=g_1\left(t, \xi ,x\left( \xi \right) ,z \right) h\left( \xi \right), \cr
		\left[ G_2\left( t,x,y,z \right) h \right] \left( \xi \right) :=g_2\left( t,\xi ,x\left( \xi \right) ,y\left( \xi \right) ,z \right) h\left( \xi \right),\nonumber
		\end{gather}
		due to {\rm (A3)} and {\rm (A4)}, for any fixed $\left( t,z\right) \in\left( \mathbb{R},\mathbb{Z}\right) ,$ the mappings
		\begin{gather}
		B_1\left( t,\cdot \right) : \mathbb{H}\times \mathbb{H}\rightarrow \mathbb{H},\quad B_2\left( t,\cdot \right) : \mathbb{H}\times \mathbb{H}\rightarrow \mathbb{H},\cr
		F_1\left( t,\cdot \right) : \mathbb{H}\rightarrow \mathcal{L}\left( \mathbb{H} \right) ,\ \ F_2\left( t,\cdot \right) : \mathbb{H}\times \mathbb{H}\rightarrow \mathcal{L}\left( \mathbb{H} \right), \cr
		G_1\left( t,\cdot ,z \right) : \mathbb{H}\rightarrow \mathcal{L}\left( \mathbb{H} \right) ,\ \ G_2\left( t,\cdot ,z \right) : \mathbb{H}\times \mathbb{H}\rightarrow \mathcal{L}\left( \mathbb{H} \right), \nonumber
		\end{gather}
		are Lipschitz continuous and linear growth conditions.}
\end{rem}

Now, for $ i=1,2 $, we define
$$
\gamma_i \left( t,s \right) :=\int_s^t{\gamma_i \left( r \right)}dr, \quad s<t,
$$
and for any $\epsilon >0$ and $\beta \ge 0$, set
$$
U_{\beta ,\epsilon,i}\left( t,s \right) =e^{\frac{1}{\epsilon}\gamma_i \left( t,s \right) A_i-\frac{\beta}{\epsilon}\left( t-s \right)}, \quad s<t.
$$
For $\epsilon =1$, we write $U_{\beta,i}\left( t,s \right)$, and 
for $\epsilon =1$ and $\beta =0$, we write $U_i\left( t,s \right) $.

Next, for any $\epsilon >0,\beta \ge 0$ and for any $u\in \mathcal{C}\big( \left[ s,t \right] ; W^{1,p}_0(\mathcal{O}) \big) , $ we define
$$
\psi _{\beta ,\epsilon,i}\left( u;s \right) \left( r \right) =\frac{1}{\epsilon}\int_s^r{U_{\beta ,\epsilon,i}\left( r,\rho \right) L_i\left( \rho \right) u\left( \rho \right)}d\rho, \quad s<r<t.
$$
In the case $\epsilon =1$, we write $\psi _{\beta, i}\left( u;s \right) \left( r \right)$, and in the case $\epsilon =1$ and $\beta =0$, we write $\psi _{i}\left( u;s \right) \left( r \right)$.
\begin{rem}\label{rem2.3}
	{\rm By using the same argument as Chapter 5 of the book \cite{pazy2012semigroups}, we can get that there exists a unique evolution system $ U\left( t,s \right) =e^{A\int_s^t{\gamma \left( r \right)}dr}\left( 0\leq s\leq t\leq T \right)  $ for $ \gamma \left( t \right) A $.
		Moreover, we also can get that $ \psi _{\beta ,\epsilon,i}\left( u;s \right) \left( t \right)  $ is the solution of
		$$
		du\left( t \right) =\frac{1}{\epsilon}\left( A_i\left( t \right) -\beta \right) u\left( t \right) dt, \quad t>s, \ u\left( s \right) =0.
		$$}
\end{rem}

\section{A priori bounds for the solution}\label{sec-3}
With all notations introduced above, system (\ref{orginal1}) can be rewritten in the following abstract form:
\begin{eqnarray}\label{orginal2}
\begin{split}
\begin{cases}
du_{\epsilon}\left( t \right) &=\left[ A_1\left( t\right) u_{\epsilon}\left( t \right) +B_1\left(t, u_{\epsilon}\left( t \right) ,v_{\epsilon}\left( t \right) \right) \right] dt+F_1\left(t, u_{\epsilon}\left( t \right) \right) d\omega ^{Q_1}\left( t \right) \\
&\quad +\int_{\mathbb{Z}}{G_1\left(t, u_{\epsilon}\left( t \right) ,z \right)}\tilde{N}_1\left( dt,dz \right), \\
dv_{\epsilon}\left( t \right) &=\frac{1}{\epsilon}\left[ \left( A_2\left( t \right) -\alpha \right) v_{\epsilon}\left( t \right) +B_2\left( t,u_{\epsilon}\left( t \right) ,v_{\epsilon}\left( t \right) \right) \right] dt \\
&\quad +\frac{1}{\sqrt{\epsilon}}F_2\left( t,u_{\epsilon}\left( t \right) ,v_{\epsilon}\left( t \right) \right) d\omega ^{Q_2}\left( t \right) \cr
&\quad +\int_{\mathbb{Z}}{G_2}\left( t,u_{\epsilon}\left( t \right) ,v_{\epsilon}\left( t \right) ,z\right) \tilde{N}_{2}^{\epsilon}\left( dt,dz \right), \\
u_{\epsilon}\left( 0 \right) &=x, \quad v_{\epsilon}\left( 0 \right) =y.
\end{cases}
\end{split}
\end{eqnarray}

According to the Remark \ref{rem2.1}, we know that the coefficients of systems (\ref{orginal2}) satisfy global Lipschitz and linear growth conditions, and the assumptions {\rm (A1)-(A4)} are uniform with respect to $t\in \mathbb{R}$. So, using the same argument as 
\cite{cerrai2003stochastic,peszat2007stochastic,pei2017two}, it is easy to prove that for any $\epsilon>0, T>0$ and $x,y\in \mathbb{H}$, there exists a unique adapted mild solution  
$\left( u_{\epsilon},v_{\epsilon} \right)$ 
to system (\ref{orginal2}) in $L^p\left( \varOmega ;D\left( \left[ 0,T \right] ;\mathbb{H}\times\mathbb{H} \right) \right)$. This means that there exist two unique adapted processes $u_{\epsilon}$ and $v_{\epsilon}$ in $L^p\left( \varOmega ;D\left( \left[ 0,T \right] ;\mathbb{H} \right) \right)$ 
such that
\begin{eqnarray}\label{or32}
\begin{split}
\begin{cases}
u_{\epsilon}\left( t \right) &=U_{1 }\left( t,0\right)x+\psi _{1}\left( u_{\epsilon};0 \right)\left( t \right)+\int_0^t{U_{1}\left( t,r \right)B_1\left(r, u_{\epsilon}\left( r \right) ,v_{\epsilon}\left( r \right) \right)}dr\cr
&\quad+\int_0^t{U_{1}\left( t,r \right)F_1\left(r, u_{\epsilon}\left( r \right) \right)}dw^{Q_1}\left( r \right) \cr
&\quad+\int_0^t{\int_{\mathbb{Z}}{U_{1}\left( t,r \right)G_1\left(r, u_{\epsilon}\left( r \right) ,z \right)}}\tilde{N}_1\left( dr,dz \right), \\
v_{\epsilon}\left( t \right) &=U_{\alpha ,\epsilon,2}\left( t,0 \right) y+\psi _{\alpha ,\epsilon,2}\left( v_{\epsilon};0 \right)\left( t \right)  \cr
&\quad+\frac{1}{\epsilon}\int_0^t{U_{\alpha ,\epsilon,2}\left( t,r \right) B_2\left( r,u_{\epsilon}\left( r \right) ,v_{\epsilon}\left( r \right) \right)}dr \cr
&\quad +\frac{1}{\sqrt{\epsilon}}\int_0^t{U_{\alpha ,\epsilon,2}\left( t,r \right) F_2\left( r,u_{\epsilon}\left( r \right) ,v_{\epsilon}\left( r \right) \right)}dw^{Q_2}\left( r \right) \cr
&\quad +\int_0^t{\int_{\mathbb{Z}}{U_{\alpha ,\epsilon,2}\left( t,r \right) G_2\left( r,u_{\epsilon}\left( r \right) ,v_{\epsilon}\left( r \right) ,z \right)}}\tilde{N}_{2}^{\epsilon}\left( dr,dz \right).
\end{cases} 
\end{split}
\end{eqnarray}
\begin{lem}\label{lem3.1}
	Under {\rm (A1)-(A4)}, for any $p\ge 1$ and $T>0$, there exists a positive constant $c_{p,T}$, such that for any $x,y\in \mathbb{H}$ and $\epsilon \in \left( 0,1 \right]$, we have   
	\begin{eqnarray}\label{en33}
	\mathbb{E}\underset{t\in \left[ 0,T \right]}{\sup}\left\| u_{\epsilon}\left( t \right) \right\| _{\mathbb{H}}^{p}\le c_{p,T}\left( 1+\left\| x \right\| _{\mathbb{H}}^{p}+\left\| y \right\| _{\mathbb{H}}^{p} \right),
	\end{eqnarray}
	\begin{eqnarray}\label{en34}
	\int_0^T{\mathbb{E}\left\| v_{\epsilon}\left( t \right) \right\| _{\mathbb{H}}^{p}}dt\le c_{p,T}\left( 1+\left\| x \right\| _{\mathbb{H}}^{p}+\left\| y \right\| _{\mathbb{H}}^{p} \right).
	\end{eqnarray}
\end{lem} 
\para{Proof:} For fixed $\epsilon \in \left( 0,1 \right]$ and $x,y\in \mathbb{H}$, for any  $t \in \left[ 0,T \right]$, we denote
$$
\varGamma _{1,\epsilon}\left( t \right) :=\int_0^t{U_{1}\left( t,r \right)F_1\left(r, u_{\epsilon}\left( r \right) \right)}dw^{Q_1}\left( r \right), 
$$
$$
\varPsi _{1,\epsilon}\left( t \right) :=\int_0^t{\int_{\mathbb{Z}}{U_{1}\left( t,r \right)G_1\left( r,u_{\epsilon}\left( r \right) ,z \right)}}\tilde{N}_1\left( dr,dz \right). 
$$
Set $\varLambda _{1,\epsilon}\left( t \right) :=u_{\epsilon}\left( t \right) -\varGamma _{1,\epsilon}\left( t \right) -\varPsi _{1,\epsilon}\left( t \right)$, we have
\begin{eqnarray}
\frac{d}{dt}\varLambda _{1,\epsilon}\left( t \right) &=&\gamma_1\left(t \right) A_1\varLambda _{1,\epsilon}\left( t \right) + L _1\left( t \right) \left( \varLambda _{1,\epsilon}\left( t \right) +\varGamma _{1,\epsilon}\left( t \right) +\varPsi _{1,\epsilon}\left( t \right) \right)\cr
&& +B_1\left( t,\varLambda _{1,\epsilon}\left( t \right) +\varGamma _{1,\epsilon}\left( t \right) +\varPsi _{1,\epsilon}\left( t \right) ,v_{\epsilon}\left( t \right) \right), \quad \varLambda _{1,\epsilon}\left( 0 \right) =x. \nonumber
\end{eqnarray}
For any $p\ge 2$, because $B_1\left( \cdot \right)$ is Lipschitz continuous, using  Young's inequality,  we have
\begin{eqnarray}\label{en333}
\frac{1}{p}\frac{d}{dt}\left\| \varLambda _{1,\epsilon}\left( t \right) \right\| _{\mathbb{H}}^{p}
&=& \left< \gamma_1\left( t\right) A_1\varLambda _{1,\epsilon}\left( t \right) ,\varLambda _{1,\epsilon}\left( t \right) \right> _{\mathbb{H}}\left\| \varLambda _{1,\epsilon}\left( t \right) \right\| _{\mathbb{H}}^{p-2}\cr
&&+\left\langle  L _1\left( t \right) \left( \varLambda _{1,\epsilon}\left( t \right) +\varGamma _{1,\epsilon}\left( t \right) +\varPsi _{1,\epsilon}\left( t \right) \right), \varLambda _{1,\epsilon}\left( t \right) \right\rangle _{\mathbb{H}}\left\| \varLambda _{1,\epsilon}\left( t \right) \right\| _{\mathbb{H}}^{p-2} \cr
&&+\left< B_1\left( t,\varLambda _{1,\epsilon}\left( t \right) +\varGamma _{1,\epsilon}\left( t \right) +\varPsi _{1,\epsilon}\left( t \right) ,v_{\epsilon}\left( t \right) \right) \right. \cr
&&\quad \left. -B_1\left(t, \varGamma _{1,\epsilon}\left( t \right) +\varPsi _{1,\epsilon}\left( t \right) ,v_{\epsilon}\left( t \right) \right) ,\varLambda _{1,\epsilon}\left( t \right) \right> _{\mathbb{H}}\left\| \varLambda _{1,\epsilon}\left( t \right) \right\| _{\mathbb{H}}^{p-2}\cr
&&+\left< B_1\left(t, \varGamma _{1,\epsilon}\left( t \right) +\varPsi _{1,\epsilon}\left( t \right) ,v_{\epsilon}\left( t \right) \right) \right. \cr
&&\quad \left. -B_1\left(t, \varPsi _{1,\epsilon}\left( t \right) ,v_{\epsilon}\left( t \right) \right) ,\varLambda _{1,\epsilon}\left( t \right) \right> _{\mathbb{H}}\left\| \varLambda _{1,\epsilon}\left( t \right) \right\| _{\mathbb{H}}^{p-2}\cr
&&+\left< B_1\left(t, \varPsi _{1,\epsilon}\left( t \right) ,v_{\epsilon}\left( t \right) \right) ,\varLambda _{1,\epsilon}\left( t \right) \right> _{\mathbb{H}}\left\| \varLambda _{1,\epsilon}\left( t \right) \right\| _{\mathbb{H}}^{p-2}\cr
&\leq& c\left\| \varLambda _{1,\epsilon}\left( t \right) \right\| _{\mathbb{H}}^{p}+c\left\| \varGamma _{1,\epsilon}\left( t \right) \right\| _{\mathbb{H}}\left\| \varLambda _{1,\epsilon}\left( t \right) \right\| _{\mathbb{H}}^{p-1} +c\left\|  \varPsi _{1,\epsilon}\left( t \right) \right\| _{\mathbb{H}}\left\| \varLambda _{1,\epsilon}\left( t \right) \right\| _{\mathbb{H}}^{p-1}\cr
&&+c\left\| B_1\left(t, \varPsi _{1,\epsilon}\left( t \right) ,v_{\epsilon}\left( t \right) \right) \right\| _{\mathbb{H}}\left\| \varLambda _{1,\epsilon}\left( t \right) \right\| _{\mathbb{H}}^{p-1}\cr
&\leq& c_p\left\| \varLambda _{1,\epsilon}\left( t \right) \right\| _{\mathbb{H}}^{p}+c_p\left( 1+\left\| \varGamma _{1,\epsilon}\left( t \right) \right\| _{\mathbb{H}}^{p}+\left\| \varPsi _{1,\epsilon}\left( t \right) \right\| _{\mathbb{H}}^{p}+\left\| v_{\epsilon}\left( t \right) \right\| _{\mathbb{H}}^{p} \right). 
\end{eqnarray}
This implies that
$$
\left\| \varLambda _{1,\epsilon}\left( t \right) \right\| _{\mathbb{H}}^{p}\le e^{c_pt}\left\| x \right\| _{\mathbb{H}}^{p}+c_p\int_0^t{e^{c_p\left( t-r \right)}\left( 1+\left\| \varGamma _{1,\epsilon}\left( r \right) \right\| _{\mathbb{H}}^{p}+\left\| \varPsi _{1,\epsilon}\left( r \right) \right\| _{\mathbb{H}}^{p}+\left\| v_{\epsilon}\left( r \right) \right\| _{\mathbb{H}}^{p} \right)}dr.
$$
According to the definition of $ \varLambda _{1,\epsilon}\left( t \right) $, for any $t\in \left[ 0,T \right]$, we have
\begin{eqnarray}
\left\| u_{\epsilon}\left( t \right) \right\| _{\mathbb{H}}^{p}
&\leq& c_{p,T}\Big( 1+\left\| x \right\| _{\mathbb{H}}^{p}+\underset{r\in \left[ 0,T \right]}{\sup}\left\| \varGamma _{1,\epsilon}\left( r \right) \right\| _{\mathbb{H}}^{p}+\underset{r\in \left[ 0,T \right]}{\sup}\left\| \varPsi _{1,\epsilon}\left( r \right) \right\| _{\mathbb{H}}^{p} \Big)\cr &&+c_{p,T}\int_0^T{\left\| v_{\epsilon}\left( r \right) \right\| _{\mathbb{H}}^{p}}dr, \nonumber
\end{eqnarray}
so
\begin{eqnarray}\label{en35}
\mathbb{E}\underset{t\in \left[ 0,T \right]}{\sup}\left\| u_{\epsilon}\left( t \right) \right\| _{\mathbb{H}}^{p}
&\leq& c_{p,T}\Big( 1+\left\| x \right\| _{\mathbb{H}}^{p}+\mathbb{E}\underset{t\in \left[ 0,T \right]}{\sup}\left\| \varGamma _{1,\epsilon}\left( t \right) \right\| _{\mathbb{H}}^{p}+\mathbb{E}\underset{t\in \left[ 0,T \right]}{\sup}\left\| \varPsi _{1,\epsilon}\left( t \right) \right\| _{\mathbb{H}}^{p} \Big)\cr
&&+c_{p,T}\int_0^T{\mathbb{E}\left\| v_{\epsilon}\left( r \right) \right\| _{\mathbb{H}}^{p}}dr.
\end{eqnarray}
According to the Lemma 4.1 in \cite{cerrai2009khasminskii} with $ \theta =0 $, it is easy to prove that
\begin{eqnarray}\label{en36}
\mathbb{E}\underset{t\in \left[ 0,T \right]}{\sup}\left\| \varGamma _{1,\epsilon}\left( t \right) \right\| _{\mathbb{H}}^{p}\le c_{p,T}\int_0^T{\left( 1+\mathbb{E}\left\| u_{\epsilon}\left( r \right) \right\| _{\mathbb{H}}^{p} \right)}dr.
\end{eqnarray}
Due to {\rm (A3)}, using Kunita's first inequality, we get 
\begin{eqnarray}\label{en108}
\mathbb{E}\underset{t\in \left[ 0,T \right]}{\sup}\left\| \varPsi _{1,\epsilon}\left( t \right) \right\| _{\mathbb{H}}^{p}
&\leq& c_p\mathbb{E}\Big( \int_0^{T}{\int_{\mathbb{Z}}{\left\| e^{\gamma_1\left( t,r \right) A_1}G_1\left( r,u_{\epsilon}\left( r \right) ,z \right) \right\| _{\mathbb{H}}^{2}}}v_1\left( dz \right) dr \Big) ^{\frac{p}{2}}\cr
&&+c_p\mathbb{E}\int_0^{T}{\int_{\mathbb{Z}}{\left\| e^{\gamma_1\left( t,r \right)A_1}G_1\left(r, u_{\epsilon}\left( r \right) ,z \right) \right\| _{\mathbb{H}}^{p}}}v_1\left( dz \right) dr \cr
&\leq& c_{p}\mathbb{E}\Big( \int_0^{T}{\left( 1+\left\| u_{\epsilon}\left( r \right) \right\| _{\mathbb{H}}^{2} \right)}dr \Big) ^{\frac{p}{2}} +c_{p} \mathbb{E}\int_0^{T}{\left( 1+\left\| u_{\epsilon}\left( r \right) \right\| _{\mathbb{H}}^{p} \right)}dr, \cr
&\leq& c_{p,T}\int_0^T{\left( 1+\mathbb{E}\left\| u_{\epsilon}\left( r \right) \right\| _{\mathbb{H}}^{p} \right)}dr,
\end{eqnarray}
so
\begin{eqnarray}\label{en37}
\mathbb{E}\underset{t\in \left[ 0,T \right]}{\sup}\left\| \varPsi _{1,\epsilon}\left( t \right) \right\| _{\mathbb{H}}^{p}\le c_{p,T}\int_0^T{\left( 1+\mathbb{E}\left\| u_{\epsilon}\left( r \right) \right\| _{\mathbb{H}}^{p} \right)}dr.
\end{eqnarray}
Substituting (\ref{en36}) and (\ref{en37}) into (\ref{en35}), we yield 
\begin{eqnarray}\label{en38}
\mathbb{E}\underset{t\in \left[ 0,T \right]}{\sup}\left\| u_{\epsilon}\left( t \right) \right\| _{\mathbb{H}}^{p}
\leq c_{p,T}\Big( 1+\left\| x \right\| _{\mathbb{H}}^{p}+\int_0^T{\mathbb{E}\left\| v_{\epsilon}\left( r \right) \right\| _{\mathbb{H}}^{p}dr} \Big) +c_{p,T}\int_0^T{\mathbb{E}\underset{\sigma\in \left[ 0,r \right]}{\sup}\left\| u_{\epsilon}\left( \sigma \right) \right\| _{\mathbb{H}}^{p}dr}. 
\end{eqnarray} 

Now, we have to estimate
$$
\int_0^T{\mathbb{E}\left\| v_{\epsilon}\left( r \right) \right\| _{\mathbb{H}}^{p}dr}.
$$
For any $ t\in \left[ 0,T \right] $, we set
$$
\varGamma _{2,\epsilon}\left( t \right) :=\frac{1}{\sqrt{\epsilon}}\int_0^t{U_{\alpha ,\epsilon,2}\left( t,r \right) F_2\left( r,u_{\epsilon}\left( r \right) ,v_{\epsilon}\left( r \right) \right) dw^{Q_2}\left( r \right)},
$$
$$
\varPsi _{2,\epsilon}\left( t \right) :=\int_0^t{\int_{\mathbb{Z}}{U_{\alpha ,\epsilon,2}\left( t,r \right) G_2\left( r,u_{\epsilon}\left( r \right) ,v_{\epsilon}\left( r \right) ,z \right)}}\tilde{N}_{2}^{\epsilon}\left( dr,dz \right).
$$
Let $\varLambda _{2,\epsilon}\left( t \right) :=v_{\epsilon}\left( t \right) -\varGamma _{2,\epsilon}\left( t \right) -\varPsi _{2,\epsilon}\left( t \right)$,  we have
\begin{eqnarray}
\frac{d}{dt}\varLambda _{2,\epsilon}\left( t \right) &=&\frac{1}{\epsilon}\left( \gamma_2 \left( t \right) A_2-\alpha \right) \varLambda _{2,\epsilon}\left( t \right) +\frac{1}{\epsilon} L _2\left( t \right) \left( \varLambda _{2,\epsilon}\left( t \right) +\varGamma _{2,\epsilon}\left( t \right) +\varPsi _{2,\epsilon}\left( t \right) \right) \cr
&&+\frac{1}{\epsilon}B_2\left( t,u_{\epsilon}\left( t \right) ,\varLambda _{2,\epsilon}\left( t \right) +\varGamma _{2,\epsilon}\left( t \right) +\varPsi _{2,\epsilon}\left( t \right) \right), \qquad\qquad  \varLambda _{2,\epsilon}\left( 0 \right) =y.\nonumber
\end{eqnarray}
For any $ p\ge2 $, because  $ \alpha >0 $ is large enough, by proceeding as in equation (\ref{en333}), we can get 
\begin{eqnarray}
\frac{1}{p}\frac{d}{dt}\left\| \varLambda _{2,\epsilon}\left( t \right) \right\| _{\mathbb{H}}^{p}
\le -\frac{\alpha}{2\epsilon}\left\| \varLambda _{2,\epsilon}\left( t \right) \right\| _{\mathbb{H}}^{p}+\frac{c_p}{\epsilon}\left( 1+\left\| u_{\epsilon}\left( t \right) \right\| _{\mathbb{H}}^{p}+\left\| \varGamma _{2,\epsilon}\left( t \right) \right\| _{\mathbb{H}}^{p}+\left\| \varPsi _{2,\epsilon}\left( t \right) \right\| _{\mathbb{H}}^{p} \right). \nonumber
\end{eqnarray}
According to the Gronwall inequality, we have
\begin{eqnarray}
\left\| \varLambda _{2,\epsilon}\left( t \right) \right\| _{\mathbb{H}}^{p}
&\le& \frac{c_p}{\epsilon}\int_0^t{e^{-\frac{\alpha p}{2\epsilon}\left( t-r \right)}\left( 1+\left\| u_{\epsilon}\left( r \right) \right\| _{\mathbb{H}}^{p}+\left\| \varGamma _{2,\epsilon}\left( r \right) \right\| _{\mathbb{H}}^{p}+\left\| \varPsi _{2,\epsilon}\left( r \right) \right\| _{\mathbb{H}}^{p} \right)}dr\cr
&&+e^{-\frac{\alpha p}{2\epsilon}t}\left\| y \right\| _{\mathbb{H}}^{p}. \nonumber
\end{eqnarray}
According to the definition of $ \varLambda _{2,\epsilon}\left( t \right) $, for any $t\in \left[ 0,T \right]$, we yield
\begin{eqnarray}
\mathbb{E}\left\| v_{\epsilon}\left( t \right) \right\| _{\mathbb{H}}^{p}
&\leq&  c_p \mathbb{E}\left\| \varGamma _{2,\epsilon}\left( t \right) \right\| _{\mathbb{H}}^{p}+c_p\mathbb{E}\left\| \varPsi _{2,\epsilon}\left( t \right) \right\| _{\mathbb{H}}^{p}+c_pe^{-\frac{\alpha p}{2\epsilon}t}\left\| y \right\| _{\mathbb{H}}^{p}\cr
&&+\frac{c_p}{\epsilon}\int_0^t{e^{-\frac{\alpha p}{2\epsilon}\left( t-r \right)}\left( 1+\mathbb{E}\left\| u_{\epsilon}\left( r \right) \right\| _{\mathbb{H}}^{p}+\mathbb{E}\left\| \varGamma _{2,\epsilon}\left( r \right) \right\| _{\mathbb{H}}^{p}+\mathbb{E}\left\| \varPsi _{2,\epsilon}\left( r \right) \right\| _{\mathbb{H}}^{p} \right)}dr.\nonumber
\end{eqnarray}
Therefore, by integrating with respect to $t$, using Young's inequality, we obtain
\begin{eqnarray}\label{en39}
\int_0^t{\mathbb{E}\left\| v_{\epsilon}\left( r \right) \right\| _{\mathbb{H}}^{p}}dr
&\leq& c_p\left( t \right) \left\| y \right\| _{\mathbb{H}}^{p}+c_p\Big( \int_0^t{\mathbb{E}\left\| \varGamma _{2,\epsilon}\left( r \right) \right\| _{\mathbb{H}}^{p}}dr+\int_0^t{\mathbb{E}\left\| \varPsi _{2,\epsilon}\left( r \right) \right\| _{\mathbb{H}}^{p}}dr \Big) \cr
&&+\frac{c_p}{\epsilon}\int_0^t{\left( 1+\mathbb{E}\left\| u_{\epsilon}\left( r \right) \right\| _{\mathbb{H}}^{p}+\mathbb{E}\left\| \varGamma _{2,\epsilon}\left( r \right) \right\| _{\mathbb{H}}^{p}+\mathbb{E}\left\| \varPsi _{2,\epsilon}\left( r \right) \right\| _{\mathbb{H}}^{p} \right)}dr\cr
&&\qquad\quad\times\int_0^t{e^{-\frac{\alpha p}{2\epsilon}r}}dr \cr
&\leq& c_p\left( t \right) \left( 1+\left\| y \right\| _{\mathbb{H}}^{p} \right)+c_p\Big( \int_0^t{\mathbb{E}\left\| u_{\epsilon}\left( r \right) \right\| _{\mathbb{H}}^{p}}dr+\int_0^t{\mathbb{E}\left\| \varGamma _{2,\epsilon}\left( r \right) \right\| _{\mathbb{H}}^{p}}dr \cr
&&+\int_0^t{\mathbb{E}\left\| \varPsi _{2,\epsilon}\left( r \right) \right\| _{\mathbb{H}}^{p}}dr \Big).
\end{eqnarray}
According to the Burkholder-Davis-Gundy inequality, 
by proceeding as Proposition 4.2 in \cite{cerrai2009khasminskii}, we can easily get
\begin{eqnarray}\label{en311}
\int_0^t{\mathbb{E}\left\| \varGamma _{2,\epsilon}\left( r \right) \right\| _{\mathbb{H}}^{p}}dr
 \leq  c_p\left( t \right) \int_0^t{\left( 1+\mathbb{E}\left\| u_{\epsilon}\left( r \right) \right\| _{\mathbb{H}}^{p}+\mathbb{E}\left\| v_{\epsilon}\left( r \right) \right\| _{\mathbb{H}}^{p} \right)}dr.
\end{eqnarray}
Concerning the stochastic term $ \varPsi _{2,\epsilon}\left( t \right) $, using  Kunita's first inequality,  we have
\begin{eqnarray}
\mathbb{E}\left\| \varPsi _{2,\epsilon}\left( t \right) \right\| _{\mathbb{H}}^{p}
&\leq& c_p\mathbb{E}\Big(\frac{1}{\epsilon} \int_0^t{\int_{\mathbb{Z}}{\big\|  e^{-\frac{\alpha}{\epsilon}\left( t-r \right)}e^{\frac{\gamma_2\left(t,r \right)  }{\epsilon}A_2}G_2\left( r,u_{\epsilon}\left( r \right) ,v_{\epsilon}\left( r \right) ,z \right) \big\|_{\mathbb{H}}^{2}}}v_2\left( dz \right) dr \Big) ^{\frac{p}{2}} \cr
&&+\frac{c_{p}}{\epsilon}\mathbb{E} \int_0^t{\int_{\mathbb{Z}}{\big\|  e^{-\frac{\alpha}{\epsilon}\left( t-r \right)}e^{\frac{\gamma_2\left(t,r \right)  }{\epsilon}A_2}G_2\left( r,u_{\epsilon}\left( r \right) ,v_{\epsilon}\left( r \right) ,z \right) \big\|_{\mathbb{H}}^{p}}}v_2\left( dz \right) dr  \cr
&\leq&\frac{c_{p}}{\epsilon^{\frac{p}{2}}}\mathbb{E} \Big( \int_0^t{ e^{-\frac{2\alpha}{\epsilon}\left( t-r \right)}
	\int_{\mathbb{Z}}{\left\| G_2\left( r,u_{\epsilon}\left( r \right) ,v_{\epsilon}\left( r \right) ,z \right) \right\| _{\mathbb{H}}^{2}}}v_2\left( dz \right) dr \Big) ^{\frac{p}{2}} \cr
&&+\frac{c_{p}}{\epsilon}\mathbb{E} \int_0^t{e^{-\frac{\alpha p}{\epsilon}\left( t-r \right)}\int_{\mathbb{Z}}{\left\| G_2\left( r,u_{\epsilon}\left( r \right) ,v_{\epsilon}\left( r \right) ,z \right) \right\| _{\mathbb{H}}^{p}}}v_2\left( dz \right) dr  \cr
&\leq&  \frac{c_p}{\epsilon ^{\frac{p}{2}}}\mathbb{E}\int_0^t{e^{-\frac{\alpha p}{2\epsilon}\left( t-r \right)}\left( 1+\lVert u_{\epsilon}\left( r \right) \rVert _{\mathbb{H}}^{p}+\lVert v_{\epsilon}\left( r \right) \rVert _{\mathbb{H}}^{p} \right)}dr\Big( \int_0^t{e^{-\frac{\alpha}{\epsilon}\frac{p}{p-2}\left( t-r \right)}}dr \Big) ^{\frac{p-2}{2}} \cr
&&+\frac{c_{p}}{\epsilon}\mathbb{E} \int_0^t{e^{-\frac{\alpha p}{\epsilon}\left( t-r \right)}\left( 1+\left\| u_{\epsilon}\left( r \right) \right\| _{\mathbb{H}}^{p}+\left\| v_{\epsilon}\left( r \right) \right\| _{\mathbb{H}}^{p} \right)}dr\cr
&\leq&\frac{c_{p}}{\epsilon} \int_0^t{e^{-\frac{\alpha p}{2\epsilon}\left( t-r \right)}\left( 1+\mathbb{E}\left\| u_{\epsilon}\left( r \right) \right\| _{\mathbb{H}}^{p}+\mathbb{E}\left\| v_{\epsilon}\left( r \right) \right\| _{\mathbb{H}}^{p} \right)}dr.\nonumber
\end{eqnarray}
By integrating with respect to $t$ both sides and using Young's inequality, we have
\begin{eqnarray}\label{en312}
\int_0^t{\mathbb{E}\left\| \varPsi _{2,\epsilon}\left( r \right) \right\| _{\mathbb{H}}^{p}}dr 
&\leq& \frac{c_{p}}{\epsilon}\int_0^t{  e^{-\frac{\alpha p}{2\epsilon}r} }dr\cdot \int_0^t{\left( 1+\mathbb{E}\left\| u_{\epsilon}\left( r \right) \right\| _{\mathbb{H}}^{p}+\mathbb{E}\left\| v_{\epsilon}\left( r \right) \right\| _{\mathbb{H}}^{p} \right)}dr\cr
&\leq& c_{p}\left( t \right) \int_0^t{\left( 1+\mathbb{E}\left\| u_{\epsilon}\left( r \right) \right\| _{\mathbb{H}}^{p}+\mathbb{E}\left\| v_{\epsilon}\left( r \right) \right\| _{\mathbb{H}}^{p} \right)}dr.
\end{eqnarray}
Substituting (\ref{en311}) and (\ref{en312}) into (\ref{en39}), we get
$$
\int_0^t{\mathbb{E}\left\| v_{\epsilon}\left( r \right) \right\| _{\mathbb{H}}^{p}}dr\le c_{p}\left( t \right) \Big( 1+\left\| y \right\| _{\mathbb{H}}^{p}+\int_0^t{ \mathbb{E}\left\| u_{\epsilon}\left( r \right) \right\| _{\mathbb{H}}^{p}   }dr\Big)+ c_{p,T} \left( t \right)\int_0^t{\mathbb{E}\left\| v_{\epsilon}\left( r \right) \right\| _{\mathbb{H}}^{p}}dr. 
$$
As $c_{p}\left( 0 \right) =0$ and $c_{p}\left( t \right) $ is a continuous increasing function, we can fix $t_0>0$, such that for any $t\le t_0$, we have $ c_{p}\left( t \right) \le {1}/{2} $, so
\begin{eqnarray}\label{en323}
\int_0^t{\mathbb{E}\left\| v_{\epsilon}\left( r \right) \right\| _{\mathbb{H}}^{p}}dr\leq c_{p}\left( t \right) \Big( 1+\left\| y \right\| _{\mathbb{H}}^{p}+\mathbb{E}\underset{r\in \left[ 0,t \right]}{\sup}\left\| u_{\epsilon}\left( r \right) \right\| _{\mathbb{H}}^{p} \Big), 
\quad t\in \left[ 0,t_0 \right].
\end{eqnarray}
Due to (\ref{en323}) and (\ref{en38}), for any $ t\in \left[ 0,t_0 \right], $ we can get
\begin{eqnarray}
\mathbb{E}\underset{r\in \left[ 0,t \right]}{\sup}\left\| u_{\epsilon}\left( r \right) \right\| _{\mathbb{H}}^{p}&\le& c_{p,T} \left( 1+\left\| x \right\| _{\mathbb{H}}^{p}+\left\| y \right\| _{\mathbb{H}}^{p} \right)+c_{p,T}\left( t \right) \mathbb{E}\underset{r\in \left[ 0,t \right]}{\sup}\left\| u_{\epsilon}\left( r \right) \right\| _{\mathbb{H}}^{p}\cr
&& +c_{p,T}\int_0^t{\mathbb{E}\underset{\sigma\in \left[ 0,r \right]}{\sup}\left\| u_{\epsilon}\left( \sigma \right) \right\| _{\mathbb{H}}^{p}dr}.\nonumber
\end{eqnarray}
Similarly, we also can fix $0<t_1\le t_0$, such that for any $t\le t_1$, we have $c_{p,T}\left( t \right) \leq {1}/{2}$, so
\begin{eqnarray} 
	\mathbb{E}\underset{r\in \left[ 0,t \right]}{\sup}\left\| u_{\epsilon}\left( r \right) \right\| _{\mathbb{H}}^{p}\le c_{p,T}  \left( 1+\left\| x \right\| _{\mathbb{H}}^{p}+\left\| y \right\| _{\mathbb{H}}^{p} \right) +c_{p,T}\int_0^t{\mathbb{E}\underset{\sigma\in \left[ 0,r \right]}{\sup}\left\| u_{\epsilon}\left( \sigma \right) \right\| _{\mathbb{H}}^{p}dr}, \quad t\in \left[ 0,t_1 \right].\nonumber
	\end{eqnarray}
	According to the Gronwall inequality, we get
\begin{eqnarray}\label{en324}
\mathbb{E}\underset{r\in \left[ 0,t \right]}{\sup}\left\| u_{\epsilon}\left( r \right) \right\| _{\mathbb{H}}^{p}\le c_{p,T}  \left( 1+\left\| x \right\| _{\mathbb{H}}^{p}+\left\| y \right\| _{\mathbb{H}}^{p} \right), \quad t\in \left[ 0,t_1 \right]. 
\end{eqnarray}
Substituting (\ref{en324}) into (\ref{en323}), it yields
\begin{eqnarray}
\int_0^{t}{\mathbb{E}\left\| v_{\epsilon}\left( r \right) \right\| _{\mathbb{H}}^{p}}dr\le c_{p,T}  \left( 1+\left\| x \right\| _{\mathbb{H}}^{p}+\left\| y \right\| _{\mathbb{H}}^{p} \right), \quad t\in \left[ 0,t_1 \right]. 
\end{eqnarray}
For any $ p\ge2 $, by repeating this in the intervals $\left[ t_1,2t_1 \right] ,  \left[ 2t_1,3t_1 \right] $ etc., we can easily get (\ref{en34}). Substituting (\ref{en34}) into (\ref{en38}) and using the Gronwall inequality again, we yield (\ref{en33}). 
Using the H\"{o}lder inequality, we can estimate (\ref{en33}) and (\ref{en34}) for $p=1$. \qed
\begin{lem}\label{lem3.2}
	Under {\rm(A1)-(A4)}, there exists $ \bar{\theta}>0 $, such that for any $T>0, p\ge 1, x\in \mathcal{D}( ( -A_1 ) ^{\theta}) $ with $\theta \in [ 0,\bar{\theta} )$ and $y\in \mathbb{H}$, there exist a positive constant $c_{p, \theta, T}>0$ such that
	\begin{eqnarray}\label{en313}
	\underset{\epsilon \in \left( 0,1 \right]}{\sup}\mathbb{E}\underset{t\in \left[ 0,T \right]}{\sup}\left\| u_{\epsilon}\left( t \right) \right\| _{\theta}^{p}\le c_{p, \theta, T}\left( 1+\left\| x \right\| _{\theta}^{p}+\left\| y \right\| _{\mathbb{H}}^{p} \right).
	\end{eqnarray}
\end{lem}
\para{Proof:} Assuming that $x\in \mathcal{D}( ( -A_1 ) ^{\theta}) (\theta \ge 0)$, for any $ t\in \left[ 0,T \right] $, we have
\begin{eqnarray}
u_{\epsilon}\left( t \right) &=&U_{1 }\left( t,0\right)x+\psi _{1}\left( u_{\epsilon};0 \right)\left( t \right)+\int_0^t{U_{1}\left( t,r \right)B_1\left(r, u_{\epsilon}\left( r \right) ,v_{\epsilon}\left( r \right) \right)}dr\cr
&&+\int_0^t{U_{1}\left( t,r \right)F_1\left(r, u_{\epsilon}\left( r \right) \right)}dw^{Q_1}\left( r \right)\cr
&& +\int_0^t{\int_{\mathbb{Z}}{U_{1}\left( t,r \right)G_1\left(r, u_{\epsilon}\left( r \right) ,z \right)}}\tilde{N}_1\left( dr,dz \right). \nonumber
\end{eqnarray}
Concerning the  second term $ \psi _{1}\left( u_{\epsilon};0 \right)\left( t \right) $, we get
\begin{eqnarray}\label{en222}
\left\| \psi _{1}\left( u_{\epsilon};0 \right)\left( t \right) \right\| _{\theta}^{p}
&\leq& c_p\Big\| \int_0^t{(-A_1)^\theta e^{\gamma_1(t,r)A_1}  L _1(r)u_\epsilon(r)}dr\Big\|_{\mathbb{H}}^{p}  \cr
&\leq& c_{p,\theta}\Big(\int_0^t{(t-r)^{-\theta }\left\| L _1(r)u_\epsilon(r) \right\|_{\mathbb{H}}} dr\Big)^p \cr
&\leq& c_{p,\theta}\underset{r \in \left[ 0,T \right]}{\sup}\left\|  u_{\epsilon}\left( r \right)\right\|_{\mathbb{H}}^{p}\Big(\int_0^t{(t-r)^{-\theta }} dr\Big)^p \cr
&\leq& c_{p, \theta, T}\left( 1+\left\| x \right\| _{\mathbb{H}}^{p}+\left\| y \right\| _{\mathbb{H}}^{p} \right).
\end{eqnarray}
For any $p\ge2$, according to the proof of Proposition 4.3 in \cite{cerrai2009khasminskii},  
and thanks to (\ref{en33}) and (\ref{en34}), it is possible to show that there exists a $\tilde{\theta}\ge 0$, such that for any $\theta \le \tilde{\theta}\land{1}/{2}$, we have 
\begin{eqnarray}\label{en314}
\mathbb{E}\underset{t\in \left[ 0,T \right]}{\sup}\Big\| \int_0^t{U_{1}\left( t,r \right)B_1\left(r, u_{\epsilon}\left( r \right) ,v_{\epsilon}\left( r \right) \right)}dr \Big\|  _{\theta}^{p}\le c_{p,\theta ,T}\left( 1+\left\| x \right\| _{\mathbb{H}}^{p}+\left\| y \right\| _{\mathbb{H}}^{p} \right).
\end{eqnarray}  
\begin{eqnarray}\label{en315}
\mathbb{E}\underset{t\in \left[ 0,T \right]}{\sup}\left\| \varGamma _{1,\epsilon}\left( t \right) \right\| _{\theta}^{p}\le c_{p, \theta, T}\left( 1+\left\| x \right\| _{\mathbb{H}}^{p}+\left\| y \right\| _{\mathbb{H}}^{p} \right).
\end{eqnarray}
Concerning the stochastic term $ \varPsi _{1,\epsilon}\left( t \right) $, using the factorization argument, we have
$$
\varPsi _{1,\epsilon}\left( t \right) =c_{\theta}\int_0^t{\left( t-r \right) ^{\theta -1}e^{\gamma_1\left( t,r \right)A_1}\phi _{\epsilon ,\theta}\left( r \right)}dr,
$$
where
$$
\phi _{\epsilon ,\theta}\left( r \right) =\int_0^r{\int_{\mathbb{Z}}{\left( r-\sigma \right) ^{-\theta}e^{\gamma_1\left( r,\sigma \right)A_1}G_1\left(\sigma, u_{\epsilon}\left( \sigma \right) ,z \right)}}\tilde{N}_1\left( d\sigma ,dz \right). 
$$
Next, for any $ p\geq2 $, let $ \hat{\theta}=( \frac{1}{4}-\frac{1}{2p} ) \land \frac{1}{2p}$, for any $ \theta \leq \hat{\theta} $, according to \rm(A3) and  \lemref{lem3.1}, using Kunita's first inequality and the H\"{o}lder inequality,  we  get  
\small\begin{eqnarray}\label{en316}
\lVert \varPsi _{1,\epsilon}\left( t \right) \rVert _{\theta}^{p}
&\leq& c_{\theta}\Big( \int_0^t{\left( t-r \right) ^{\theta -1}\lVert \phi _{\epsilon ,\theta}\left( r \right) \rVert _{\theta}}dr\ \Big) ^p \cr
&\leq& c_{\theta}\underset{r\in \left[ 0,t \right]}{\sup}\lVert \phi _{\epsilon ,\theta}\left( r \right) \rVert _{\theta}^{p}\Big( \int_0^t{\left( t-r \right) ^{\theta -1}}dr\ \Big) ^p \cr
&\leq& c_{p,\theta ,T}\underset{r\in \left[ 0,t \right]}{\sup}\Big\| \int_0^r{\int_{\mathbb{Z}}{\left( r-\sigma \right) ^{-\theta}\left( -A_1 \right) ^{\theta}e^{\gamma_1\left( r,\sigma \right)A_1}G_1\left(\sigma, u_{\epsilon}\left( \sigma \right) ,z \right)}}\tilde{N}_1\left( d\sigma ,dz \right) \Big\| _{\mathbb{H}}^{p} \cr
&\leq& c_{p,\theta ,T}\underset{r\in \left[ 0,t \right]}{\sup}\Big( \int_0^r{\int_{\mathbb{Z}}{\left( r-\sigma \right) ^{-2\theta} \left( r-\sigma \right)^{-2\theta} \lVert G_1\left(\sigma, u_{\epsilon}\left( \sigma \right) ,z \right) \rVert _{\mathbb{H}}^{2}}}v_1\left( dz \right) d\sigma \Big) ^{\frac{p}{2}}\cr
&&+c_{p,\theta ,T}\underset{r\in \left[ 0,t \right]}{\sup}\int_0^r{\int_{\mathbb{Z}}{\left( r-\sigma \right) ^{-p\theta} \left( r-\sigma \right)^{-p\theta}\lVert G_1\left(\sigma, u_{\epsilon}\left( \sigma \right) ,z \right) \rVert _{\mathbb{H}}^{p}}}v_1\left( dz \right) d\sigma \cr
&\leq& c_{p,\theta ,T}\underset{r\in \left[ 0,t \right]}{\sup}\bigg[    \int_{0}^{r}\Big(\int_{\mathbb{Z}}{\lVert G_1\left(\sigma, u_{\epsilon}\left( \sigma \right) ,z \right) \rVert _{\mathbb{H}}^{2}v_1\left( dz \right)} \Big)^\frac{p}{2}d\sigma \cr
&&\qquad\qquad\qquad\qquad\qquad\qquad\qquad\qquad\qquad\quad\times\Big( \int_0^r{\left( r-\sigma \right) ^{-\frac{4 p\theta}{p-2}}d\sigma} \Big) ^{\frac{p-2}{2}}\bigg] \cr
&&+c_{p,\theta ,T}\underset{r\in \left[ 0,t \right]}{\sup}\bigg[ \int_0^r{\left( r-\sigma \right) ^{-2p\theta}}d\sigma\cdot  \Big( \underset{\sigma \in \left[ 0,r \right]}{\sup}\int_{\mathbb{Z}}{\lVert G_1\left(\sigma, u_{\epsilon}\left( \sigma \right) ,z \right) \rVert _{\mathbb{H}}^{p}v_1\left( dz \right)} \Big)\bigg]  \cr
&\leq&c_{p,\theta ,T}\Big( 1+\underset{\sigma \in \left[ 0,T \right]}{\sup}\left\|  u_{\epsilon}\left( \sigma \right)\right\|_{\mathbb{H}}^{p} \Big).\nonumber
\end{eqnarray}
So, due to (\ref{en33}), we yield
\begin{eqnarray}
\mathbb{E}\underset{t\in \left[ 0,T \right]}{\sup}\lVert \varPsi _{1,\epsilon}\left( t \right) \rVert _{\theta}^{p}\leq c_{p, \theta, T}\left( 1+\left\| x \right\| _{\mathbb{H}}^{p}+\left\| y \right\| _{\mathbb{H}}^{p} \right). 
\end{eqnarray}
Hence, if we choose $\bar{\theta}:={1}/{8}\land \tilde{\theta}\land \hat{\theta}$, thanks to (\ref{en222}), (\ref{en314}), (\ref{en315}) and (\ref{en316}), for any $p\ge 2$ and $\theta <\bar{\theta}$, we get
\begin{eqnarray}
\mathbb{E}\underset{t\in \left[ 0,T \right]}{\sup}\left\| u_{\epsilon}\left( t \right) \right\| _{\theta}^{p}
&\leq& c_{p, \theta, T}\left( 1+\left\| x \right\| _{\theta}^{p}+\left\| y \right\| _{\mathbb{H}}^{p} \right). \nonumber
\end{eqnarray}
Using the H\"{o}lder inequality, we can estimate (\ref{en313}) for $p=1$. \qed
\begin{lem}\label{lem3.3}
	Under {\rm (A1)-(A4)}, for any $\theta \in [0,\bar{\theta}) $ and $  0\le h\le 1,$ there exists $\beta \left( \theta \right) >0,$ such that, for any $T>0,p\ge 1,x\in \mathcal{D}( ( -A_1 ) ^{\theta}),y\in \mathbb{H}$ and $ t\in \left[ 0,T \right] $, it holds
	\begin{eqnarray}\label{en317}
	\underset{\epsilon \in \left( 0,1 \right]}{\sup}\mathbb{E}\left\| u_{\epsilon}\left( t \right) -u_{\epsilon}\left( t+h \right) \right\| _{\mathbb{H}}^{p}\leq c_{p, \theta, T}\big( h^{\beta \left( \theta \right) p}+h\big) \left( 1+\left\| x \right\| _{\theta}^{p}+\left\| y \right\| _{\mathbb{H}}^{p} \right). 
	\end{eqnarray}
\end{lem}
\para{Proof:} For any $t\ge 0 $ with $t,t+h\in \left[ 0,T \right] $, we have
\begin{eqnarray}\label{en318}
u_{\epsilon}\left( t+h \right) -u_{\epsilon}\left( t \right) &=&\left(U_1(t+h,t)-I \right) u_{\epsilon}\left( t \right) +\psi _{1}\left( u_{\epsilon};t \right)\left( t+h \right)\cr
&&+\int_t^{t+h}{U_1(t+h,r) B_1\left(r, u_{\epsilon}\left( r \right) ,v_{\epsilon}\left( r \right) \right) }dr\cr
&&+\int_t^{t+h}{U_1(t+h,r)}F_1\left(r, u_{\epsilon}\left( r \right) \right)dw^{Q_1}\left( r \right) \cr
&&+\int_t^{t+h}{\int_{\mathbb{Z}}{U_1(t+h,r)}G_1\left(r, u_{\epsilon}\left( r \right) ,z \right)}\tilde{N}_1\left( dr,dz \right)\cr
&:=&\sum_{i=1}^5{\mathcal{I}_{t}^{i}}. 
\end{eqnarray}
By proceeding as the proof of Proposition 4.4 in \cite{cerrai2009khasminskii} and (\ref{en222}), fix $\theta \in[0,\bar{\theta} ) $, for any $p\ge 1$, it is possible to show that
\begin{eqnarray}\label{en319}
\mathbb{E}\left\| \mathcal{I}_{t}^{1} \right\| _{\mathbb{H}}^{p}
\leq c_{p, \theta, T}h^{\theta p}\left( 1+\left\| x \right\| _{\theta}^{p}+\left\| y \right\| _{\mathbb{H}}^{p} \right),
\end{eqnarray}
\begin{eqnarray}\label{en320}
\mathbb{E}\left\| \mathcal{I}_{t}^{2} \right\| _{\mathbb{H}}^{p}
\leq c_{p,T}h^{p-1}\left( 1+\left\| x \right\| _{\mathbb{H}}^{p}+\left\| y \right\| _{\mathbb{H}}^{p} \right), 
\end{eqnarray}
\begin{eqnarray}\label{en223}
\mathbb{E}\left\| \mathcal{I}_{t}^{3} \right\| _{\mathbb{H}}^{p}
\leq c_{p, T}h^p\left( 1+\left\| x \right\| _{\mathbb{H}}^{p}+\left\| y \right\| _{\mathbb{H}}^{p} \right),
\end{eqnarray}
\begin{eqnarray}\label{en321}
\mathbb{E}\left\| \mathcal{I}_{t}^{4} \right\| _{\mathbb{H}}^{p}
\leq c_{p,T}h^{\frac{p-2}{2}-\frac{\beta _1\left( \rho _1-2 \right)}{\rho _1}\frac{p}{2}}\left( 1+\left\| x \right\| _{\mathbb{H}}^{p}+\left\| y \right\| _{\mathbb{H}}^{p} \right). 
\end{eqnarray}	
According to the proof of (\ref{en108}), using the H\"{o}lder inequality and (\ref{en33}), we have
\begin{eqnarray}\label{en322}
\mathbb{E}\left\| \mathcal{I}_{t}^{5} \right\| _{\mathbb{H}}^{p}
&\leq& c_{p}\mathbb{E}\Big( \int_t^{t+h}{\left( 1+\left\| u_{\epsilon}\left( r \right) \right\| _{\mathbb{H}}^{2} \right)}dr \Big) ^{\frac{p}{2}} +c_{p} \mathbb{E} \int_t^{t+h}{\left( 1+\left\| u_{\epsilon}\left( r \right) \right\| _{\mathbb{H}}^{p} \right)}dr, \cr
&\leq& c_{p} \big( h^{\frac{p-2}{2}}+1\big)  \int_t^{t+h}{\left( 1+\mathbb{E}\left\| u_{\epsilon}\left( r \right) \right\| _{\mathbb{H}}^{p} \right)}dr  \cr
&\leq& c_{p,T}\big( h^{\frac{p}{2}}+h\big) \left( 1+\left\| x \right\| _{\mathbb{H}}^{p}+\left\| y \right\| _{\mathbb{H}}^{p} \right).  
\end{eqnarray}
Then, if we take $\bar{p}>1$, such that
$$
\frac{\beta _1\left( \rho _1-2 \right)}{\rho _1}\frac{\bar{p}}{\bar{p}-2}<1,
$$
we can get
\begin{eqnarray}
\mathbb{E}\left\| u_{\epsilon}\left( t+h \right) -u_{\epsilon}\left( t \right) \right\| _{\mathbb{H}}^{p}
&\leq& c_{p,T}\Big( h^{(1-\frac{1}{\bar{p}})p}+h^p+h^{\left( \frac{1}{2}-\frac{1}{\bar{p}}-\frac{\beta _1\left( \rho _1-2 \right)}{2\rho _1} \right) p}+h^{\frac{p}{2}}+h  \Big) \cr
&& \times\left( 1+\left\| x \right\| _{\mathbb{H}}^{p}+\left\| y \right\| _{\mathbb{H}}^{p} \right) 
+c_{p,\theta,T}h^{\theta p}\left( 1+\left\| x \right\| _{\theta}^{p}+\left\| y \right\| _{\mathbb{H}}^{p} \right). \nonumber
\end{eqnarray}
As we are assuming $\left|  h \right| \le 1$, (\ref{en317}) follows for any $p\ge \bar{p}$ by taking
$$
\beta \left( \theta \right) :=\min \left\{ \theta, 1-\frac{1}{\bar{p}},1, \frac{1}{2}-\frac{1}{\bar{p}}-\frac{\beta _1\left( \rho _1-2 \right)}{2\rho _1},\frac{1}{2} \right\}. 
$$
From the H\"{o}lder inequality, we can estimate (\ref{en317}) for $p<\bar{p}$, 
$$
\mathbb{E}\left\| u_{\epsilon}\left( t+h \right) -u_{\epsilon}\left( t \right) \right\| _{\mathbb{H}}^{p}\le \left[ \mathbb{E}\left\| u_{\epsilon}\left( t+h \right) -u_{\epsilon}\left( t \right) \right\| _{\mathbb{H}}^{\bar{p}} \right] ^{\frac{p}{\bar{p}}}.
$$
So, we have (\ref{en317}).\qed

According to the above lemma, we get that for every $\epsilon \in \left( 0,1 \right] , $ the function $ u_{\epsilon}\left( t \right) $ is uniformly bounded about $ t\in\left[ 0,T\right] , $ and it is also equicontinuous at every point of $ t\in \left[0,T \right]  $.  In view of the Theorem 12.3 in \cite{Billingsley1968Convergence}, we can infer that the set $\left\{   u_{\epsilon}  \right\} _{\epsilon \in \left( 0,1 \right]}$ is relatively compact in $ D\left(\left[0,T \right]; \mathbb{H}\right). $ In addition, according to the Theorem 13.2 in \cite{Billingsley1968Convergence}  and the above lemma, by using the Chebyshev's inequality, we also can get that the family of probability measures $ \left\lbrace \mathcal{L} \left( u_\epsilon \right)  \right\rbrace_{\epsilon\in \left( 0,1\right] }  $ is tight in $ \mathcal{P} \left( D\left(\left[0,T \right]; \mathbb{H}\right)\right) . $

\section{An evolution family of measures for the fast equation}\label{sec-4}
For any frozen slow component $x\in \mathbb{H}$, any initial condition $y\in \mathbb{H}$, and any $s\in \mathbb{R}$, we introduce the following problem
\begin{eqnarray}\label{en41}
dv\left( t \right) &=&\left[ \left( A_2\left( t \right) -\alpha \right) v\left( t \right) +B_2\left( t,x,v\left( t \right) \right) \right] dt+F_2\left( t,x,v\left( t \right) \right) d\bar{\omega}^{Q_2}\left( t \right) \cr
&&+{\int_{\mathbb{Z}}{G_2}\left( t,x,v\left( t \right) ,z \right)}{\tilde{N}_{{2}^{'}}}\left( dt,dz \right), \qquad\qquad\qquad\qquad v\left( s \right)=y,
\end{eqnarray}
where 
$$
\bar{w}^{Q_2}\left( t \right) =\left\{ \begin{array}{l}
w_{1}^{Q_2}\left( t \right),  \\
w_{2}^{Q_2}\left( -t \right),  \\
\end{array} \right. \begin{array}{c}
if\ t\ge 0,\\
if\ t<0,\\
\end{array}
$$
$$
{\tilde{N}_{{2}^{'}}}\left( t,z \right) =\left\{ \begin{array}{l}
{\tilde{N}_{{1}^{'}}}\left( t,z \right), \\
{\tilde{N}_{{3}^{'}}}\left( -t,z \right),  \\
\end{array} \right. \begin{array}{c}
if\ t\ge 0,\\
if\ t<0,\\
\end{array}
$$
for two independent $Q_2$-Wiener processes $w_{1}^{Q_2}\left( t \right)$, $w_{2}^{Q_2}\left( t \right)$ and two independent compensated Poisson measures ${\tilde{N}_{{1}^{'}}}\left( dt,dz \right)$, ${\tilde{N}_{{3}^{'}}}\left( dt,dz \right)$ with the same  L\'{e}vy measure are both defined as in Section \ref{sec-2}. 

According to the definition of the operator $ \psi _{\alpha,2}\left( \cdot ;s \right)  $, we know that the mapping $ \psi _{\alpha,2}\left( \cdot ;s \right) :\mathcal{C}\left( \left[ s,T \right] ;\mathbb{H} \right) $ $\rightarrow \mathcal{C}\left( \left[ s,T \right] ;\mathbb{H} \right) $ is a linear bounded operator and it is Lipschitz continuous. Hence, we have that, for any $ x,y\in \mathbb{H}, p\ge 1$ and $ s<T $, there exists a unique mild solution  $ v^x\left( \cdot ;s,y \right)\in L^p\left( \varOmega ;D\left( \left[ 0,T \right] ;\mathbb{H} \right) \right)   $ (\cite{peszat2007stochastic})  in the following form 
\begin{eqnarray}
v^x\left( t;s,y \right) &=&U_{\alpha,2}\left( t,s \right) y+\psi _{\alpha,2}\left( v^x\left( \cdot;s,y \right) ;s \right)\left( t \right) +\int_s^t{U_{\alpha,2}\left( t,r \right) B_2\left( r,x,v^x\left( r;s,y \right) \right)}dr\cr
&&+\int_s^t{U_{\alpha,2}\left( t,r \right) F_2\left( r,x,v^x\left( r;s,y \right) \right)}d\bar{w}^{Q_2}\left( r \right) \cr
&&+\int_s^t{\int_{\mathbb{Z}}{U_{\alpha,2}\left( t,r \right) G_2\left( r,x,v^x\left( r;s,y \right) ,z \right)}}\tilde{N}_{{2}^{'}}\left( dr,dz \right).\nonumber
\end{eqnarray}

Moreover, if the space $ D\left( \mathbb{R};\mathbb{H} \right) $ endowed with the topology of uniform convergence on bounded intervals, an $ \left\{ \mathcal{F}_t \right\} _{t\in \mathbb{R}} $-adapted process $ v^x \in L^p\left( \varOmega ;D\left( \left[ 0,T \right] ;\mathbb{H} \right) \right)$ is a mild solution of the equation
\begin{eqnarray}\label{en42}
dv\left( t \right) &=&\left[ \left( A_2\left( t \right) -\alpha \right) v\left( t \right) +B_2\left( t,x,v\left( t \right) \right) \right] dt+F_2\left( t,x,v\left( t \right) \right) d\bar{\omega}^{Q_2}\left( t \right) \cr
&&+\int_{\mathbb{Z}}{G_2}\left( t,x,v\left( t \right) ,z \right) \tilde{N}_{{2}^{'}}\left( dt,dz \right),
\end{eqnarray}
where $ t\in \mathbb{R} $. Then, for every $ s<t $, we have
\begin{eqnarray}
v^x\left( t \right) &=&U_{\alpha,2}\left( t,s \right) v^x\left( s \right)+\psi _{\alpha,2}\left( v^x ;s \right)\left( t \right) +\int_s^t{U_{\alpha,2}\left( t,r \right) B_2\left( r,x,v^x\left( r \right) \right)}dr\cr
&&+\int_s^t{U_{\alpha,2}\left( t,r \right) F_2\left( r,x,v^x\left( r \right) \right)}d\bar{w}^{Q_2}\left( r \right) \cr
&&+\int_s^t{\int_{\mathbb{Z}}{U_{\alpha,2}\left( t,r \right) G_2\left( r,x,v^x\left( r \right) ,z \right)}}\tilde{N}_{{2}^{'}}\left( dr,dz \right). \nonumber
\end{eqnarray}

In what follows, for any $ x\in \mathbb{H} $ and any adapted process $v $, we set
\begin{eqnarray}\label{en43}
\varGamma _{\alpha}\left( v;s \right) \left( t \right) :=\int_s^t{U_{\alpha,2}\left( t,r \right) F_2\left( r,x,v\left( r \right) \right)}d\bar{w}^{Q_2}\left( r \right), \quad t>s,
\end{eqnarray}
\begin{eqnarray}\label{en44}
\varPsi _{\alpha}\left( v;s \right) \left( t \right) :=\int_s^t{\int_{\mathbb{Z}}{U_{\alpha,2}\left( t,r \right) G_2\left( r,x,v\left( r \right) ,z \right)}}\tilde{N}_{{2}^{'}}\left( dr,dz \right), \quad t>s.
\end{eqnarray}
For any $   0<\delta <\alpha  $ and any $  v_1,v_2 $  with $ s<t $, by proceeding as in the proof of Lemma 7.1 in \cite{cerrai2006asymptotic}, it is possible to show that there exists $ \bar{p}>1 $, such that for any $ p\ge \bar{p},$  we have
\begin{eqnarray}\label{en45}
\underset{r\in \left[ s,t \right]}{\sup}e^{\delta p\left( r-s \right)}\mathbb{E}\left\| \varGamma _{\alpha}\left( v_1;s \right) \left( r \right) -\varGamma _{\alpha}\left( v_2;s \right) \left( r \right) \right\| _{\mathbb{H}}^{p}
&\leq& c_{p,1}\frac{L_{f_2}^{p}}{\left( \alpha -\delta \right) ^{c_{p,2}}} \cr
&&\times\underset{r\in \left[ s,t \right]}{\sup}e^{\delta p\left( r-s \right)}\mathbb{E}\left\| v_1\left( r \right) -v_2\left( r \right) \right\| _{\mathbb{H}}^{p}, \cr
&&
\end{eqnarray}	
where $  L_{f_2}   $ is the Lipschitz constant of $ f_2 $, and $ c_{p,1},c_{p,2} $ are two suitable positive constants independent of $ \alpha >0 $ and $ s<t $. 

For the stochastic term  $\varPsi _{\alpha}\left( v;s \right) \left( t \right)$, using Kunita's first inequality \cite[Theorem 4.4.23]{Applebaum2009Processes}, we get
\begin{small}
	\begin{eqnarray}
	&&\mathbb{E}\left\| \varPsi _{\alpha}\left( v_1;s \right) \left( t \right) -\varPsi _{\alpha}\left( v_2;s \right) \left( t \right) \right\| _{\mathbb{H}}^{p}\cr
	&\leq&c_p \mathbb{E}\Big(\int_s^t{\int_{\mathbb{Z}}{\left\|e^{-\alpha \left( t-r \right)} e^{\gamma_2 \left( t,r \right)A_2}\left[ G_2\left( r,x,v_1\left( r \right) ,z \right) -G_2\left( r,x,v_2\left( r \right) ,z \right) \right] \right\| _{\mathbb{H}}^{2}}}v_{2^{'}}\left(dz \right)  dr\Big)^{\frac{p}{2}}\cr
	&&+c_p \mathbb{E} \int_s^t{\int_{\mathbb{Z}}{\left\| e^{- \alpha  \left( t-r \right)}e^{\gamma_2 \left( t,r \right)A_2}\left[ G_2\left( r,x,v_1\left( r \right) ,z \right) -G_2\left( r,x,v_2\left( r \right) ,z \right) \right] \right\| _{\mathbb{H}}^{p}}}v_{2^{'}}\left(dz \right)   dr\cr
	&\leq& c_p\mathbb{E}\Big( \int_s^t{\int_{\mathbb{Z}}{e^{-2\alpha \left( t-r \right)}\left\| G_2\left( r,x,v_1\left( r \right) ,z \right) -G_2\left( r,x,v_2\left( r \right) ,z \right) \right\| _{\mathbb{H}}^{2}}}v_{2^{{'}}}\left( dz \right) dr\Big)^{\frac{p}{2}}\cr
	&&+c_p\mathbb{E}\int_s^t{\int_{\mathbb{Z}}{e^{- \alpha p \left( t-r \right)}\left\| G_2\left( r,x,v_1\left( r \right) ,z \right) -G_2\left( r,x,v_2\left( r \right) ,z \right) \right\| _{\mathbb{H}}^{p}}}v_{2^{{'}}}\left( dz \right)dr\cr
	&\leq& c_pL_{g_2}^{p}\Big( \int_s^t{e^{-2\left( \alpha -\delta \right) \left( t-r \right)}e^{-2\delta \left( t-s \right)}e^{2\delta \left( r-s \right)}\mathbb{E}\left\| v_1\left( r \right) -v_2\left( r \right) \right\| _{\mathbb{H}}^{2}}dr \Big)^{\frac{p}{2}}\cr
	&&+ c_pL_{g_2}^{p} \int_s^t{e^{-p\left( \alpha -\delta \right) \left( t-r \right)}e^{-\delta p\left( t-s \right)}e^{\delta p \left( r-s \right)}\mathbb{E}\left\| v_1\left( r \right) -v_2\left( r \right) \right\| _{\mathbb{H}}^{p}}dr \cr
	&\leq& c_pL_{g_2}^{p}\Big[ \Big( \int_s^t{e^{-2\left( \alpha -\delta \right) \left( t-r \right)}}dr \Big)^{\frac{p}{2}} + \int_s^t{e^{-p\left( \alpha -\delta \right) \left( t-r \right)}}dr    \Big]\cr
	&&\qquad\qquad\qquad\qquad\qquad\qquad\qquad\qquad\quad \times e^{-\delta p\left( t-s \right)}\underset{r\in \left[ s,t \right]}{\sup}e^{\delta p\left( r-s \right)}\mathbb{E}\left\| v_1\left( r \right) -v_2\left( r \right) \right\| _{\mathbb{H}}^{p}\cr
	&\leq& c_{p,1}\frac{L_{g_2}^{p}}{\left( \alpha -\delta \right) ^{c_{p,2}}}e^{-\delta p\left( t-s \right)}\underset{r\in \left[ s,t \right]}{\sup}e^{\delta p\left( r-s \right)}\mathbb{E}\left\| v_1\left( r \right) -v_2\left( r \right) \right\| _{\mathbb{H}}^{p},\nonumber
	\end{eqnarray}
\end{small}
so
\begin{eqnarray}\label{en46}
\underset{r\in \left[ s,t \right]}{\sup}e^{\delta p\left( r-s \right)}\mathbb{E}\left\| \varPsi _{\alpha}\left( v_1;s \right) \left( r \right) -\varPsi _{\alpha}\left( v_2;s \right) \left( r \right) \right\| _{\mathbb{H}}^{p}
&\leq& c_{p,1}\frac{L_{g_2}^{p}}{\left( \alpha -\delta \right) ^{c_{p,2}}}\cr
&&\times\underset{r\in \left[ s,t \right]}{\sup}e^{\delta p\left( r-s \right)}\mathbb{E}\left\| v_1\left( r \right) -v_2\left( r \right) \right\| _{\mathbb{H}}^{p},\cr
&&
\end{eqnarray}	
where  $ L_{g_2} $ is the Lipschitz constant of $ g_2 $, and $ c_{p,1},c_{p,2} $  are two suitable positive constants independent of $ \alpha >0 $  and $ s<t $.

Moreover, using {\rm (A4)}, we can show that
\begin{align}\label{en47}
\underset{r\in \left[ s,t \right]}{\sup}e^{\delta p\left( r-s \right)}\mathbb{E}\left\| \varGamma _{\alpha}\left( v;s \right) \left( r \right) \right\|_{\mathbb{H}}^{p} \le c_{p,1}\frac{M_{f_2}^{p}}{\left( \alpha -\delta \right) ^{c_{p,2}}}\underset{r\in \left[ s,t \right]}{\sup}e^{\delta p\left( r-s \right)}\left( 1+\mathbb{E}\left\| v\left( r \right) \right\| _{\mathbb{H}}^{p} \right),
\end{align}
\begin{align}\label{en48}
\underset{r\in \left[ s,t \right]}{\sup}e^{\delta p\left( r-s \right)}\mathbb{E}\left\| \varPsi _{\alpha}\left( v;s \right) \left( r \right) \right\|_{\mathbb{H}}^{p} \le c_{p,1}\frac{M_{g_2}^{p}}{\left( \alpha -\delta \right) ^{c_{p,2}}}\underset{r\in \left[ s,t \right]}{\sup}e^{\delta p\left( r-s \right)}\left( 1+\mathbb{E}\left\| v\left( r \right) \right\| _{\mathbb{H}}^{p} \right),
\end{align}
where $ M_{f_2}, M_{g_2} $ are the linear growth constants of $ f_2 , g_2 $, and $ c_{p,1},c_{p,2} $  are two suitable positive constants independent of $ \alpha >0 $  and $ s<t $.

For any fixed adapted process $ v $, let us introduce the problem
\begin{eqnarray}\label{en49}
d\rho \left( t \right) &=&\left( A_2\left( t \right) -\alpha \right) \rho \left( t \right) dt+F_2\left( t,x,v \left( t \right) \right) d\bar{\omega}^{Q_2}\left( t \right) \cr
&&+\int_{\mathbb{Z}}{G_2}\left( t,x,v \left( t \right) ,z \right) \tilde{N}_{{2}^{'}}\left( dt,dz \right), \qquad \rho \left( s \right) =0.
\end{eqnarray}
We denote that its unique mild solution is $ \rho _{\alpha}\left( v;s \right) $. This means that $ \rho _{\alpha}\left( v;s \right)  $  solves the equation
$$
\rho _{\alpha}\left( v;s \right) \left( t \right) =\psi _{\alpha,2}\left( \rho _{\alpha}\left( v;s \right) ;s \right) \left( t \right) +\varGamma _{\alpha}\left( v;s \right) \left( t \right) +\varPsi _{\alpha}\left( v;s \right) \left( t \right), \quad s<t<T.
$$
Due to (\ref{en45}) and (\ref{en46}), using the same arguement as the equation (5.8) in \cite{cerrai2017averaging}, it is easy to prove that for any process $  v_1, v_2 $ and $ 0<\delta<\alpha $, we have
\begin{eqnarray}\label{en411}
\underset{r\in \left[ s,t \right]}{\sup}e^{\delta p\left( r-s \right)}\mathbb{E}\left\| \rho _{\alpha}\left( v_1;s \right) \left( r \right) -\rho _{\alpha}\left( v_2;s \right) \left( r \right) \right\| _{\mathbb{H}}^{p}
&\leq& c_{p,1}\frac{L^p}{\left( \alpha -\delta \right) ^{c_{p,2}}}\cr
&&\times\underset{r\in \left[ s,t \right]}{\sup}e^{\delta p\left( r-s \right)}\mathbb{E}\left\| v_1\left( r \right) -v_2\left( r \right) \right\| _{\mathbb{H}}^{p},\cr
&&
\end{eqnarray}
where $ L=\max \left\{ L_{b_2},L_{f_2},L_{g_2} \right\} $.

Similarly, thanks to (\ref{en47}) and (\ref{en48}), for any process $  v  $ and $ 0<\delta<\alpha $, we  can prove that 
\begin{align}\label{en412}
\underset{r\in \left[ s,t \right]}{\sup}e^{\delta p\left( r-s \right)}\mathbb{E}\left\| \rho _{\alpha}\left( v;s \right) \left( r \right) \right\| _{\mathbb{H}}^{p}\le c_{p,1}\frac{M^p}{\left( \alpha -\delta \right) ^{c_{p,2}}}\underset{r\in \left[ s,t \right]}{\sup}e^{\delta p\left( r-s \right)}\mathbb{E}\left( 1+\left\| v\left( r \right) \right\| _{\mathbb{H}}^{p} \right),
\end{align}
where $ M=\max \left\{ M_{b_2},M_{f_2},M_{g_2} \right\} $.
\begin{lem}\label{lem4.1}
	Under {\rm (A1)-(A4)}, there exists $ \delta >0 $,  such that for any $ x,y\in \mathbb{H} $ and $ p\ge 1 $,
	\begin{eqnarray}\label{en413}
	\mathbb{E}\left\| v^x\left( t;s,y \right) \right\| _{\mathbb{H}}^{p}\le c_p\left( 1+\left\| x \right\| _{\mathbb{H}}^{p}+e^{-\delta p\left( t-s \right)}\left\| y \right\| _{\mathbb{H}}^{p} \right), \quad s<t.
	\end{eqnarray}	
\end{lem}
\para{Proof:} We set $ \varLambda _{\alpha}\left( t \right) :=v^x\left( t;s,y \right) -\rho _{\alpha}\left( t \right) $, where $ \rho _{\alpha}\left( t \right) =\rho _{\alpha}\left( v^x\left( \cdot ;s,y \right) ;s \right) \left( t \right) $  is the solution of the problem (\ref{en49}) with $ v=v^x\left( \cdot ;s,y \right) $. Using Young's inequality, we have
\begin{eqnarray}
\frac{1}{p}\frac{d}{dt}\left\| \varLambda _{\alpha}\left( t \right) \right\| _{\mathbb{H}}^{p}
&\leq& \left< \left( A_2\left( t \right)-\alpha \right) \varLambda _{\alpha}\left( t \right) ,\varLambda _{\alpha}\left( t \right) \right> _{\mathbb{H}}\left\| \varLambda _{\alpha}\left( t \right) \right\| _{\mathbb{H}}^{p-2}\cr
&&+\left< B_2\left( t,x,\varLambda _{\alpha}\left( t \right) +\rho _{\alpha}\left( t \right) \right) -B_2\left( t,x,\rho _{\alpha}\left( t \right) \right) ,\varLambda _{\alpha}\left( t \right) \right> _{\mathbb{H}}\left\| \varLambda _{\alpha}\left( t \right) \right\| _{\mathbb{H}}^{p-2}\cr
&&+\left< B_2\left( t,x,\rho _{\alpha}\left( t \right) \right) ,\varLambda _{\alpha}\left( t \right) \right> _{\mathbb{H}}\left\| \varLambda _{\alpha}\left( t \right) \right\| _{\mathbb{H}}^{p-2}\cr
&\leq&-\alpha \left\| \varLambda _{\alpha}\left( t \right) \right\| _{\mathbb{H}}^{p}+c\left\| \varLambda _{\alpha}\left( t \right) \right\| _{\mathbb{H}}^{p}+c\left( 1+\left\| x \right\| _{\mathbb{H}}+\left\| \rho _{\alpha}\left( t \right) \right\| _{\mathbb{H}} \right)\left\| \varLambda _{\alpha}\left( t \right) \right\| _{\mathbb{H}}^{p-1}\cr
&\leq& -\alpha \left\| \varLambda _{\alpha}\left( t \right) \right\| _{\mathbb{H}}^{p}+c_p\left\| \varLambda _{\alpha}\left( t \right) \right\| _{\mathbb{H}}^{p}+c_p\left( 1+\left\| x \right\| _{\mathbb{H}}^{p}+\left\| \rho _{\alpha}\left( t \right) \right\| _{\mathbb{H}}^{p} \right). \nonumber
\end{eqnarray}
Because $ \alpha $ is large enough,  we can find $ \eta =\alpha -c_p >0 $, such that
$$
\frac{d}{dt}\left\| \varLambda _{\alpha}\left( t \right) \right\| _{\mathbb{H}}^{p}\le -\eta p\left\| \varLambda _{\alpha}\left( t \right) \right\| _{\mathbb{H}}^{p}+c_p\left( 1+\left\| x \right\| _{\mathbb{H}}^{p}+\left\| \rho _{\alpha}\left( t \right) \right\| _{\mathbb{H}}^{p} \right).
$$
According to the Gronwall inequality, we have
$$
\left\| \varLambda _{\alpha}\left( t \right) \right\| _{\mathbb{H}}^{p}\le e^{-\eta p\left( t-s \right)}\left\| y \right\| _{\mathbb{H}}^{p}+c_p\left( 1+\left\| x \right\| _{\mathbb{H}}^{p} \right) +c_p\int_s^t{e^{-\eta p\left( t-r \right)}\left\| \rho _{\alpha}\left( r \right) \right\| _{\mathbb{H}}^{p}}dr.
$$
So, for any  $ p\ge 1 $,
\begin{eqnarray}
\left\| v^x\left( t;s,y \right) \right\| _{\mathbb{H}}^{p}&\leq& c_p\left\| \rho _{\alpha}\left( t \right) \right\| _{\mathbb{H}}^{p}+ c_p e^{-\eta p\left( t-s \right)}\left\| y \right\| _{\mathbb{H}}^{p}+ c_p \left( 1+\left\| x \right\| _{\mathbb{H}}^{p} \right) \cr
&&+ c_p \int_s^t{e^{-\eta p\left( t-r \right)}\left\| \rho _{\alpha}\left( r \right) \right\| _{\mathbb{H}}^{p}}dr.\nonumber
\end{eqnarray}
Fix  $ 0<\delta <\eta $, according to (\ref{en412}), we get
\begin{eqnarray}
e^{\delta p\left( t-s \right)}\mathbb{E}\left\| v^x\left( t;s,y \right) \right\| _{\mathbb{H}}^{p}
&\leq& c_pe^{\delta p\left( t-s \right)}\mathbb{E}\left\| \rho _{\alpha}\left( t \right) \right\| _{\mathbb{H}}^{p}+c_pe^{p\left( \delta -\eta \right) \left( t-s \right)}\left\| y \right\| _{\mathbb{H}}^{p}\cr
&&+c_pe^{\delta p\left( t-s \right)}\left( 1+\left\| x \right\| _{\mathbb{H}}^{p} \right) +c_p\int_s^t{e^{\delta p\left( r-s \right)}\mathbb{E}\left\| \rho _{\alpha}\left( r \right) \right\| _{\mathbb{H}}^{p}}dr\cr
&\leq& c_{p,1}\frac{M^p}{\left( \alpha -\delta \right) ^{c_{p,2}}}\underset{r\in \left[ s,t \right]}{\sup}e^{\delta p\left( r-s \right)}\left( 1+\mathbb{E}\left\| v^x\left( r;s,y \right) \right\| _{\mathbb{H}}^{p} \right) \cr
&&+c_pe^{p\left( \delta -\eta \right) \left( t-s \right)}\left\| y \right\| _{\mathbb{H}}^{p}+c_pe^{\delta p\left( t-s \right)}\left( 1+\left\| x \right\| _{\mathbb{H}}^{p} \right) \cr
&&+c_{p,1}\frac{M^p}{\left( \alpha -\delta \right) ^{c_{p,2}}}\int_s^t{\underset{r \in \left[ s,t \right]}{\sup}e^{\delta p\left( r-s \right)}\left( 1+\mathbb{E}\left\| v^x\left( r;s,y \right) \right\| _{\mathbb{H}}^{p} \right)}dr.\nonumber
\end{eqnarray}
Taking  $ \alpha _1=\left( 2c_{p,1}M^p \right) ^{\frac{1}{c_{p,2}}}+\delta $, when  $ \alpha \ge \alpha _1 ,$ we have 
\begin{eqnarray}
\underset{r\in \left[ s,t \right]}{\sup}e^{\delta p\left( r-s \right)}\mathbb{E}\left\| v^x\left( r;s,y \right) \right\| _{\mathbb{H}}^{p}
&\leq& c_p\left\| y \right\| _{\mathbb{H}}^{p}+c_pe^{\delta p\left( t-s \right)}\left( 1+\left\| x \right\| _{\mathbb{H}}^{p} \right) \cr
&&+\int_s^t{\underset{r \in \left[ s,t \right]}{\sup}e^{\delta p\left( r-s \right)}\mathbb{E}\left\| v^x\left( r;s,y \right) \right\| _{\mathbb{H}}^{p}}dr.\nonumber
\end{eqnarray}
Due to the Gronwall lemma, we have
\begin{eqnarray}
\underset{r\in \left[ s,t \right]}{\sup}e^{\delta p\left( r-s \right)}\mathbb{E}\left\| v^x\left( r;s,y \right) \right\| _{\mathbb{H}}^{p}
\le c_p\left\| y \right\| _{\mathbb{H}}^{p}+c_pe^{\delta p\left( t-s \right)}\left( 1+\left\| x \right\| _{\mathbb{H}}^{p} \right). \nonumber
\end{eqnarray}
Hence, we get (\ref{en413}).\qed

\begin{lem}\label{lem4.2}
	Under {\rm (A1)-(A4)}, 
	for any $ t\in \mathbb{R} $  and  $ x,y\in \mathbb{H} $, for all $ p\ge 1 $, there exists $ \eta ^x\left( t \right) \in L^p\left( \varOmega ;\mathbb{H} \right) $ such that
	\begin{eqnarray}\label{en414}
	\underset{s\rightarrow -\infty}{\lim}\mathbb{E}\left\| v^x\left( t;s,y \right) -\eta ^x\left( t \right) \right\| _{\mathbb{H}}^{p}=0.
	\end{eqnarray}
	Moreover, for any  $ p\ge 1 $, there exists some  $ \delta _p>0 $, such that
	\begin{eqnarray}\label{en415}
	\mathbb{E}\left\| v^x\left( t;s,y \right) -\eta ^x\left( t \right) \right\| _{\mathbb{H}}^{p}\le c_pe^{-\delta _p\left( t-s \right)}\left( 1+\left\| x \right\| _{\mathbb{H}}^{p}+\left\| y \right\| _{\mathbb{H}}^{p} \right). 
	\end{eqnarray}
	Finally,  $ \eta ^x $ is a mild solution in $ \mathbb{R} $ of equation (\ref{en42}).
\end{lem}
\para{Proof:} Fix  $ h>0 $ and define
$$
\rho \left( t \right) =v^x\left( t;s,y \right) -v^x\left( t;s-h,y \right), \quad s<t.
$$
We know that  $ \rho \left( t \right) $ is the unique mild solution of the problem
\begin{eqnarray}\label{en416}
\begin{split}
\begin{cases}
d\rho \left( t \right) &=\left[ \left( A_2\left( t \right) -\alpha \right) \rho \left( t \right) +B_2\left( t,x,v^x\left( t;s,y \right) \right) -B_2\left( t,x,v^x\left( t;s-h,y \right) \right) \right] dt\\
&\quad+\left[ F_2\left( t,x,v^x\left( t;s,y \right) \right) -F_2\left( t,x,v^x\left( t;s-h,y \right) \right) \right] d\bar{w}^{Q_2}\left( t \right) \\
&\quad+\int_{\mathbb{Z}}{\left[ G_2\left( t,x,v^x\left( t;s,y \right) \right) -G_2\left( t,x,v^x\left( t;s-h,y \right) \right) \right]}\tilde{N}_{{2}^{'}}\left( dt,dz \right)\\
\rho\left(s \right)&=y-v^x\left( s;s-h,y \right), 
\end{cases} 
\end{split}
\end{eqnarray}
and 
\begin{eqnarray}
\rho\left(t \right) &=&U_{\alpha,2} \left( t,s \right) \left(y-v^x\left( s;s-h,y \right)\right) +\psi _{\alpha,2}\left( \rho;s \right) \left( t \right) \cr
&&+\int_s^t{U_{\alpha,2}\left( t,r \right) \left[B_2\left( r,x,v^x\left( r;s,y \right) \right) -B_2\left( r,x,v^x\left( r;s-h,y \right) \right) \right]}dr \cr
&&+\int_s^t{U_{\alpha,2}\left( t,r \right) \left[F_2\left( r,x,v^x\left( r;s,y \right) \right) -F_2\left( r,x,v^x\left( r;s-h,y \right) \right) \right]}d\bar{w}^{Q_2}\left( r \right) \cr
&&+\int_s^t{\int_{\mathbb{Z}}{U_{\alpha,2}\left( t,r \right) \left[G_2\left( r,x,v^x\left( r;s,y \right) \right) -G_2\left( r,x,v^x\left( r;s-h,y \right) \right) \right]}}\tilde{N}_{{2}^{'}}\left( dr,dz \right). \nonumber
\end{eqnarray}
Multiply both sides of the above equation by $ e^{\delta p\left( t-s \right)} $. Because $ \alpha $ large enough, according to the Lemma 2.4 in \cite{cerrai2017averaging}, we have
\begin{eqnarray}
e^{\delta p\left( t-s \right)}\mathbb{E} \left\| \rho\left(t \right) \right\| _{\mathbb{H}}^{p}
&\leq& c_pe^{ \left( \delta-\alpha\right)   p\left( t-s \right)}\mathbb{E}\left\| e^{\gamma_2 \left( t,s \right)A_2} \left(y-v^x\left( s;s-h,y \right)\right) \right\| _{\mathbb{H}}^{p}\cr
&&+c_p \mathbb{E}\Big\|  \int_s^t{e^{\gamma_2 \left( t,r \right)A_2}e^{\left( \delta -\alpha \right) \left( t-r \right)}e^{\delta \left( r-s \right)}\left[B_2\left( r,x,v^x\left( r;s,y \right) \right)\right.  }\cr
&&\qquad\qquad\qquad\qquad {\left. -B_2\left( r,x,v^x\left( r;s-h,y \right) \right) \right]}dr \Big\| _{\mathbb{H}}^{p}\cr
&&+c_p\mathbb{E}\Big\|  \int_s^t{e^{\gamma_2 \left( t,r \right)A_2}e^{\left( \delta -\alpha \right) \left( t-r \right)}e^{\delta \left( r-s \right)}\left[F_2\left( r,x,v^x\left( r;s,y \right) \right)\right. }\cr
&&\qquad\qquad\qquad\qquad {\left.  -F_2\left( r,x,v^x\left( r;s-h,y \right) \right) \right]}d\bar{w}^{Q_2}\left( r \right) \Big\| _{\mathbb{H}}^{p}\cr
&&+c_p \mathbb{E}\Big\| \int_s^t{\int_{\mathbb{Z}}{e^{\gamma_2 \left( t,r \right)A_2}e^{\left( \delta -\alpha \right) \left( t-r \right)}e^{\delta \left( r-s \right)}\left[G_2\left( r,x,v^x\left( r;s,y \right) \right) \right. }}\cr
&&\qquad\qquad\qquad\qquad {{\left. -G_2\left( r,x,v^x\left( r;s-h,y \right) \right) \right]}}\tilde{N}_{{2}^{'}}\left( dr,dz \right) \Big\|  _{\mathbb{H}}^{p}\cr
&:=&\sum_{i=1}^4{\mathcal{I}_{t}^{i}}.\nonumber
\end{eqnarray}
According to the Lemma 3.1 in \cite{cerrai2009khasminskii}, we know that for any $J\in \mathcal{L}\left( L^{\infty}\left( D \right) ,\mathbb{H} \right) \cap \mathcal{L}\left( \mathbb{H},L^1\left( D \right) \right)$ with $ J=J^{\ast} $ and $ s\ge0 $, we have
\begin{eqnarray}\label{en310}
\left\| e^{sA_i}JQ_i \right\| _{2}^{2}\le K_is^{-\frac{\beta _i \left( \rho _i-2 \right)}{\rho _i}}e^{-\frac{\alpha \left( \rho _i+2 \right)}{\rho _i}s}\left\| J \right\| _{\mathcal{L}\left( L^{\infty}\left( D \right) ,\mathbb{H} \right)}^{2},
\end{eqnarray}
where
$$
K_i=\big( {\beta _i}/{e} \big) ^{\frac{\beta _i\left( \rho _i-2 \right)}{\rho _i}}\zeta _{i}^{\frac{\left( \rho _i-2 \right)}{\rho _i}}\kappa _{i}^{\frac{2}{\rho _i}}.
$$
Taking  $ \bar{p}>1 $, such that  $ \frac{\beta _2\left( \rho _2-2 \right)}{\rho _2}\frac{\bar{p}}{\bar{p}-2}<1 $. Then, using the Burkholder-Davis-Gundy inequality and Kunita's first inequality,  we can get that for any  $p\ge \bar{p}$ and $ 0<\delta<\alpha $, it yields
\begin{eqnarray}
\mathcal{I}_{t}^{2}
&\leq& c_pL_{b_2}^{p}\underset{r\in \left[ s,t \right]}{\sup}e^{\delta p \left( r-s \right)}\mathbb{E}\left\| \rho\left( r\right) \right\| _{\mathbb{H}}^{p}\cdot
\Big( \int_s^t{ e^{\left( \delta -\alpha \right) \left( t-r \right)}  }dr \Big)^{p}\cr
&\leq& c_p\frac{L_{b_2}^{p}}{\left( \alpha -\delta\right)^p}\underset{r\in \left[ s,t \right]}{\sup}e^{\delta p \left( r-s \right)}\mathbb{E}\left\| \rho\left( r\right) \right\| _{\mathbb{H}}^{p},\nonumber
\end{eqnarray}
\small\begin{eqnarray}
\mathcal{I}_{t}^{3}
&\leq& c_p 
\mathbb{E}\Big( \int_s^t{\left\| e^{\gamma_2 \left( t,r \right)A_2}e^{\left( \delta -\alpha \right) \left( t-r \right)}e^{\delta \left( r-s \right)}\left[F_2\left( r,x,v^x\left( r;s,y \right) \right)
	\right. \right. }\cr
&&\qquad\qquad\qquad\qquad\qquad\qquad\qquad\qquad\qquad\qquad \left. \left. -F_2\left( r,x,v^x\left( r;s-h,y \right) \right) \right]Q_2 \right\|_{2}^{2}dr \Big)^{\frac{p}{2}}\cr
&\leq&c_pL_{f_2}^{p}K_{2}^{\frac{p}{2}}\Big( \int_s^t{e^{2\left( \delta -\alpha \right) \left( t-r \right)} \gamma_2 \left( t,r \right) ^{-\frac{\beta _2\left( \rho _2-2 \right)}{\rho _2}}e^{-\frac{\alpha \left( \rho _2+2 \right)}{\rho _2} \gamma_2 \left( t,r \right)}}dr \Big)^{\frac{p}{2}}\cr
&&\qquad\qquad\qquad\qquad\qquad\qquad\qquad\qquad\qquad\qquad\qquad\times\underset{r\in \left[ s,t \right]}{\sup}e^{\delta p \left( r-s \right)}\mathbb{E}\left\|\rho\left( r\right) \right\| _{\mathbb{H}}^{p} \cr
&\leq&c_pL_{f_2}^{p}\underset{r\in \left[ s,t \right]}{\sup}e^{\delta p \left( r-s \right)}\mathbb{E}\left\| \rho\left( r\right) \right\| _{\mathbb{H}}^{p}
\cdot\Big( \int_s^{t}{e^{2\left( \delta -\alpha \right) \left( t-r \right)}\left[ \gamma_0\left( t-r \right) \right] ^{-\frac{\beta _2\left( \rho _2-2 \right)}{\rho _2}}}dr \Big)^{\frac{p}{2}}\cr
&\leq&c_pL_{f_2}^{p}\underset{r\in \left[ s,t \right]}{\sup}e^{\delta p \left( r-s \right)}\mathbb{E}\left\| \rho\left( r\right) \right\| _{\mathbb{H}}^{p} \cdot\Big( \int_0^{t-s}{r^{-\frac{\beta _2\left( \rho _2-2 \right)}{\rho _2}\frac{p}{p-2}}}dr \Big)^{\frac{p-2}{2}}\Big( \int_0^{t-s}{e^{-p\left( \alpha -\delta \right) r}}dr \Big)\cr
&\leq& c_p\frac{L_{f_2}^{p}}{\alpha -\delta}\underset{r\in \left[ s,t \right]}{\sup}e^{\delta p \left( r-s \right)}\mathbb{E}\left\| \rho\left( r\right) \right\| _{\mathbb{H}}^{p},\nonumber
\end{eqnarray}
\begin{eqnarray}
\mathcal{I}_{t}^{4}
&\leq& c_p\mathbb{E} \Big( \int_s^t{\int_{\mathbb{Z}}{\big\| e^{\gamma_2 \left( t,r \right)A_2}e^{\left( \delta -\alpha \right) \left( t-r \right)}e^{\delta \left( r-s \right)}\big[  G_2\left( r,x,v^x\left( r;s,y \right) \right)}} \cr
&&\qquad\qquad\qquad\qquad\qquad\qquad{{-G_2\left( r,x,v^x\left( r;s-h,y \right) \right)  \big]\big\| _{\mathbb{H}}^{2}}}v_{{2}^{'}}\left( dz \right) dr \Big)^{\frac{p}{2}}\cr
&&+c_p\mathbb{E}\int_s^t{\int_{\mathbb{Z}}{\big\|   e^{\gamma_2 \left( t,r \right)A_2}e^{\left( \delta -\alpha \right) \left( t-r \right)}e^{\delta \left( r-s \right)}\big[  G_2\left( r,x,v^x\left( r;s,y \right) \right)}} \cr
&&\qquad\qquad\qquad\qquad\qquad\qquad{{-G_2\left( r,x,v^x\left( r;s-h,y \right) \right)  \big]  \big\|_{\mathbb{H}}^{p}}}v_{{2}^{'}}\left( dz \right) dr \cr
&\leq& c_pL_{g_2}^{p}\underset{r\in \left[ s,t \right]}{\sup}e^{\delta p \left( r-s \right)}\mathbb{E}\left\| \rho\left( r\right) \right\| _{\mathbb{H}}^{p}
\cdot\Big[ \Big( \int_0^{t-s}{e^{-2\left( \alpha -\delta \right) r}}dr \Big)^{\frac{p}{2}}+\int_0^{t-s}{e^{-p\left( \alpha -\delta \right) r}} dr\Big] \cr
&\leq& c_p\frac{L_{g_2}^{p}}{\left( \alpha -\delta \right) ^{c_p}}\underset{r\in \left[ s,t \right]}{\sup}e^{\delta p \left( r-s \right)}\mathbb{E}\left\| \rho\left( r\right) \right\| _{\mathbb{H}}^{p}.\nonumber
\end{eqnarray}
Hence, we have
$$
\underset{r\in \left[ s,t \right]}{\sup}e^{\delta p\left( r-s \right)}\mathbb{E}\left\| \rho\left( r\right) \right\| _{\mathbb{H}}^{p}\le c_p\left\|y-v^x\left( s;s-h,y \right) \right\| _{\mathbb{H}}^{p}+c_{p,1}\frac{L^p}{\left( \alpha -\delta \right) ^{c_{p,2}}}\underset{r\in \left[ s,t \right]}{\sup}e^{\delta p\left( r-s \right)}\mathbb{E}\left\| \rho\left( r\right) \right\| _{\mathbb{H}}^{p}.
$$
Therefore, for $ \alpha >0 $ large enough, we can find  $ 0<\bar{\delta}_p<\alpha  $, such that
$$
c_{p,1}\frac{L}{\left( \alpha -\bar{\delta}_p \right) ^{c_{p,2}}}<1.
$$
This implies that
$$
\underset{r\in \left[ s,t \right]}{\sup}e^{p\bar{\delta}_p\left( r-s \right)}\mathbb{E}\left\| \rho\left( r\right) \right\| _{\mathbb{H}}^{p}\le c_p\left\|y-v^x\left( s;s-h,y \right) \right\| _{\mathbb{H}}^{p}.
$$
Let $ \delta _p=p\bar{\delta}_p $, thanks to the \lemref{lem4.1}, we have
\begin{eqnarray}\label{en419}
\mathbb{E}\left\| v^x\left( t;s,y \right) -v^x\left( t;s-h,y \right) \right\| _{\mathbb{H}}^{p} 
&\le& c_pe^{-\delta _p\left( t-s \right)}\left\| y-v^x\left( s;s-h,y \right) \right\| _{\mathbb{H}}^{p}\cr
&\le& c_pe^{-\delta _p\left( t-s \right)}\left( 1+\left\| x \right\| _{\mathbb{H}}^{p}+\left\| y \right\| _{\mathbb{H}}^{p}+e^{-\delta ph}\left\| y \right\| _{\mathbb{H}}^{p} \right).\cr
&&
\end{eqnarray}
Because $ L^p\left( \varOmega ;\mathbb{H} \right) $ is completeness, for any $ p\geq\bar{p} $, let  $ s\rightarrow -\infty $ in (\ref{en419}),  there exists $ \eta ^x\left( t \right) \in L^p\left( \varOmega ;\mathbb{H} \right) $  such that (\ref{en414}) hold. Then, if we let  $ h\rightarrow \infty  $ in (\ref{en419}), we obtain (\ref{en415}). Using the H\"{o}lder inequality, we can get (\ref{en414}) and (\ref{en415}) holds for any $p<\bar{p}$.

If we take  $ y_1,y_2\in \mathbb{H} $, use the same arguments for $ v^x\left( t;s,y_1 \right) -v^x\left( t;s,y_2 \right), s<t $,  we have
$$
\mathbb{E}\left\| v^x\left( t;s,y_1 \right) -v^x\left( t;s,y_2 \right) \right\| _{\mathbb{H}}^{p}\le c_pe^{-\delta _p\left( t-s \right)}\left\| y_1-y_2 \right\| _{\mathbb{H}}^{p}, \quad s<t.
$$
Let  $ s\rightarrow -\infty  $, this means that the limit $ \eta ^x\left( t \right) $  does not depend on the initial condition  $ y\in \mathbb{H} $.

Finally, we prove that $ \eta ^x\left( t \right) $  is a mild solution of equation (\ref{en42}). Due to the limit $ \eta ^x\left( t \right) $  does not depend on the initial condition, we can let initial condition $ y=0 $. For any $ s<t $  and  $ h>0 $,  we have 
\begin{eqnarray}
v^x\left( t;s-h,0 \right)
&=&U_{\alpha,2}\left( t,s \right) v^x\left( s;s-h,0 \right) +\psi _{\alpha,2}\left( v^x\left( \cdot ;s-h,0 \right) ;s \right) \left( t \right) \cr
&&+\int_s^t{U_{\alpha,2}\left( t,r \right) B_2\left( r,x,v^x\left( r;s-h,0 \right) \right)}dr\cr
&&+\int_s^t{U_{\alpha,2}\left( t,r \right) F_2\left( r,x,v^x\left( r;s-h,0 \right) \right)}d\bar{w}^{Q_2}\left( r \right) \cr
&&+\int_s^t{\int_{\mathbb{Z}}{U_{\alpha,2}\left( t,r \right) G_2\left( r,x,v^x\left( r;s-h,0 \right) ,z \right)}}\tilde{N}_{{2}^{'}}\left( dr,dz \right).\nonumber
\end{eqnarray} 
Let $ h\rightarrow \infty  $ on both sides, due to (\ref{en414}), we can get, for any  $ s<t $, have
\begin{eqnarray}\label{en424}
\eta ^x\left( t \right) &=&U_{\alpha,2}\left( t,s \right) \eta ^x\left( s \right) +\psi _{\alpha,2}\left( \eta ^x;s \right) \left( t \right) +\int_s^t{U_{\alpha,2}\left( t,r \right) B_2\left( r,x,\eta ^x\left( r \right) \right)}dr\cr
&&+\int_s^t{U_{\alpha,2}\left( t,r \right) F_2\left( r,x,\eta ^x\left( r \right) \right)}d\bar{w}^{Q_2}\left( r \right) \cr
&&+\int_s^t{\int_{\mathbb{Z}}{U_{\alpha,2}\left( t,r \right) G_2\left( r,x,\eta ^x\left( r \right) ,z \right)}}\tilde{N}_{{2}^{'}}\left( dr,dz \right),
\end{eqnarray}
this means that $ \eta ^x\left( t \right) $  is a mild solution of equation (\ref{en42}).\qed

For any  $ t\in \mathbb{R} $ and  $ x\in \mathbb{H} $, we denote that the law of the random variable  $ \eta ^x\left( t \right)  $ is $ \mu _{t}^{x} $, and we introduce the transition evolution operator
$$
P_{s,t}^{x}\varphi \left( y \right) =\mathbb{E}\varphi \left( v^x\left( t;s,y \right) \right), \quad s<t, \ y\in \mathbb{H},
$$
where $ \varphi \in \mathcal{B}_b\left( \mathbb{H} \right) $.

Due to (\ref{en413}) and (\ref{en414}), for any  $ p\ge 1 $, we have
\begin{eqnarray}\label{en426}
\underset{t\in \mathbb{R}}{\sup}\mathbb{E}\left\| \eta ^x\left( t \right) \right\| _{\mathbb{H}}^{p}\le c_p\left( 1+\left\| x \right\| _{\mathbb{H}}^{p} \right),\quad x\in \mathbb{H},
\end{eqnarray}
so that
\begin{eqnarray}\label{en427}
\underset{t\in \mathbb{R}}{\sup}\int_{\mathbb{H}}{\left\| y \right\| _{\mathbb{H}}^{p}}\mu _{t}^{x}\left( dy \right) \le c_p\left( 1+\left\| x \right\| _{\mathbb{H}}^{p} \right),\quad x\in \mathbb{H}.
\end{eqnarray}

According to the above conclusion, by using the same arguments as \cite[Proposition 5.3]{cerrai2017averaging}, we know that the family $\left\{ \mu _{t}^{x} \right\}_{t \in \mathbb{R} }$ defines an evolution system of probability measures on $ \mathbb{H} $ for equation (\ref{en41}). This means that for any $ t\in \mathbb{R} $, $ \mu _{t}^{x} $ is a probability measure on $ \mathbb{H} $, and it holds that
\begin{eqnarray}\label{en428}
\int_{\mathbb{H}}P_{s,t}^{x}\varphi \left( y \right)\mu _{s}^{x}\left( dy \right) =\int_{\mathbb{H}}{\varphi \left( y \right)}\mu _{t}^{x}\left( dy \right), \quad s<t,
\end{eqnarray}
for every $ \varphi \in \mathcal{C}_b(\mathbb{H}) $. Moreover, we also have 
\begin{eqnarray}\label{en429}
\Big| P_{s,t}^{x}\varphi \left( y \right) -\int_{\mathbb{H}}{\varphi \left( y \right)}\mu _{t}^{x}\left( dy \right) \Big|\le ce^{-\delta_1 \left( t-s \right)}\left( 1+\lVert x \rVert _{\mathbb{H}} \right).
\end{eqnarray}

\begin{lem}
	Under {\rm (A1)-(A4)}, the family of measures 
	\begin{eqnarray}\label{0109}
	\varLambda :=\left\{ \mu _{t}^{x}:t\in \mathbb{R},\ x\in \mathbb{H} \right\}  
	\end{eqnarray}
	is tight in $ \mathcal{P} \left( D\left(\left[0,T \right]; \mathbb{H}\right)\right) . $
\end{lem}
\para{Proof:} 
%
According to the define of $ v^x\left( \cdot;s,0 \right), $ for any $ t>s$ and  $h>0$, we have
\begin{eqnarray}\label{en0318}
v^x\left( t+h;s,0 \right) -v^x\left( t;s,0 \right) &=&\left(U_{\alpha,2}(t+h,t)-I \right) v^x\left( t;s,0 \right) +\psi _{\alpha,2}\left( v^x\left( \cdot;s,0 \right);t \right)\left( t+h \right)\cr
&&+\int_t^{t+h}{U_{\alpha,2}(t+h,r) B_2\left(r, x, v^x\left( r;s,0 \right) \right) }dr\cr
&&+\int_t^{t+h}{U_{\alpha,2}(t+h,r)}F_2\left(r, x, v^x\left( r;s,0 \right) \right)d\bar{w}^{Q_2}\left( r \right) \cr
&&+\int_t^{t+h}{\int_{\mathbb{Z}}{U_{\alpha,2}(t+h,r)}G_2\left(r, x, v^x\left( r;s,0 \right) ,z \right)}\tilde{N}_{{2}^{'}}\left( dr,dz \right)\nonumber
\end{eqnarray}
Due to the assumption {\rm (A4)}, we know that the mappings $ B_2,  F_2, G_2 $ are linearly growing. For any $ t>s,\ 0<h<1 $ and $ x\in \mathbb{H}, $ analogous to the proof of \lemref{lem3.3} and using the estimate (\ref{en413}), we can get 
\begin{eqnarray}\label{en0317}
\mathbb{E}\left\| v^x\left( t+h;s,0 \right) -v^x\left( t;s,0 \right) \right\| _{\mathbb{H}}^{p}\leq c_{p}\big( h^{\kappa\left( p \right) }+h\big) \left( 1+\left\| x \right\|_{\mathbb{H}}^{p}  \right), 
\end{eqnarray}
where $ \kappa\left( p \right) $ is a function of $ p $ and satisfies $ \kappa\left( p \right)>0. $

From the \lemref{lem4.2}, we know that the limit $ \eta ^x\left( t \right) $  does not depend on the initial condition, we can get
\begin{eqnarray}
\mathbb{E}\lVert \eta ^x\left( t+h \right) -\eta ^x\left( t \right) \rVert _{\mathbb{H}}^{p}&\le& \mathbb{E}\lVert \eta ^x\left( t+h \right) -v^x\left( t+h;s,0 \right) \rVert _{\mathbb{H}}^{p}+\mathbb{E}\lVert v^x\left( t+h;s,0 \right) -v^x\left( t;s,0 \right) \rVert _{\mathbb{H}}^{p}\cr
&&+\mathbb{E}\lVert v^x\left( t;s,0 \right) -\eta ^x\left( t \right) \rVert _{\mathbb{H}}^{p}.\nonumber
\end{eqnarray}
Then, thanks to (\ref{en415}) and (\ref{en0317}), let $ s\rightarrow \infty, $   we have
\begin{eqnarray}\label{en431}
\mathbb{E}\left\| \eta ^{x}\left(t+h \right) -\eta ^{x}\left(t \right)  \right\| _{\mathbb{H}}^{p}\le c_{p}\big( h^{\kappa\left( p \right) }+h\big) \left( 1+\left\| x \right\|_{\mathbb{H}}^{p}  \right). 
\end{eqnarray}
Moreover, we have proved that the family $\left\{ \mu _{t}^{x} \right\}_{t \in \mathbb{R} }$ defines an evolution system of probability measures for equation (\ref{en41}). By using the Chebyshev's inequality and equation (\ref{en426}), it yields 
\begin{eqnarray}\label{en051}
\underset{n\rightarrow \infty}{\lim}\mu _{t}^{x}\left( \eta _x\left( t \right) :\lVert \eta _x\left( t \right) \rVert _{\mathbb{H}}\geq n \right)\le \underset{n\rightarrow \infty}{\lim} \frac{\mathbb{E}\left\| \eta ^x\left( t \right) \right\|_{\mathbb{H}}^{2}}{n^2} \le \underset{n\rightarrow \infty}{\lim}c_2\frac{1+\lVert x \rVert _{\mathbb{H}}^{2}}{n^2}=0.
\end{eqnarray}
In addition, for any $ t>s,\ h>0 $ and fixed $ \epsilon>0 $, according to the equation (\ref{en431}), we have
\begin{eqnarray}\label{en052}
\underset{h\rightarrow 0}{\lim}\mu _{t}^{x}\left( \eta _x\left( t \right) :\lVert \eta _x\left( t+h \right) -\eta _x\left(t \right) \rVert _{\mathbb{H}}\geq \epsilon \right) &\le& \underset{h\rightarrow 0}{\lim} \frac{\mathbb{E}\left\| \eta ^{x}\left(t+h \right) -\eta ^{x}\left(t \right)  \right\| _{\mathbb{H}}^{2}}{\epsilon ^2}\cr
&\le& \underset{h\rightarrow 0}{\lim}c_{2 }\frac{\big( h^{ \kappa \left( 2 \right)}+h \big) \big( 1+\lVert x \rVert _{\mathbb{H}}^{2} \big)}{\epsilon ^2}=0.
\end{eqnarray}	
According to the Theorem 13.2 in \cite{Billingsley1968Convergence}, 
the equations (\ref{en051})  and (\ref{en052}) imply that 
the family of measures 
\begin{eqnarray}\label{109}
\varLambda :=\left\{ \mu _{t}^{x}:t\in \mathbb{R},x\in \mathbb{H} \right\}  
\end{eqnarray}
is tight in $ \mathcal{P} \left( D\left(\left[0,T \right]; \mathbb{H}\right)\right) . $\qed 

In order to get the averaged equation, we must ensure the existence of the averaged coefficient $ \bar{B}_1 $. So, we need the evolution family of measures satisfying some nice properties. We give the following assumption.
\begin{enumerate}[({A}5)]
	\item 
	\begin{enumerate}
		\item The functions $ \gamma_2:\mathbb{R}\rightarrow \left( 0,\infty \right)  $  and $ l_2:\mathbb{R}\times\mathcal{O}\rightarrow \mathbb{R}^d $  are periodic, with the same period.
		\item The families of functions
		\begin{gather}
		\mathbf{B}_{1,R}:=\left\{ b_1\left( \cdot ,\xi ,\sigma \right) :\ \xi \in\mathcal{O},\ \sigma \in B_{\mathbb{R}^2}\left( R \right) \right\},\cr
		\mathbf{B}_{2,R}:=\left\{ b_2\left( \cdot ,\xi ,\sigma \right) :\ \xi \in\mathcal{O},\ \sigma \in B_{\mathbb{R}^2}\left( R \right) \right\},\cr
		\mathbf{F}_R:=\left\{ f_2\left( \cdot ,\xi ,\sigma \right) :\ \xi \in\mathcal{O},\ \sigma \in B_{\mathbb{R}^2}\left( R \right) \right\},\cr
		\mathbf{G}_R:=\left\{ g_2\left( \cdot ,\xi ,\sigma ,z \right) :\ \xi \in\mathcal{O},\ \sigma \in B_{\mathbb{R}^2}\left( R \right) ,\ z\in \mathbb{Z} \right\},\nonumber
		\end{gather}
		are uniformly almost periodic for any  $ R>0 $. 
	\end{enumerate}
\end{enumerate}
\begin{rem}\label{rem4.5}
	{\rm Similar with the proof of Lemma 6.2 in \cite{cerrai2017averaging}, we get that under {\rm (A5)}, for any  $ R>0 $, the families of functions
		\begin{gather}
		\left\{ B_1\left( \cdot ,x,y \right) :\left( x,y \right) \in B_{\mathbb{H}\times \mathbb{H}}\left( R \right) \right\}, \quad
		\left\{ B_2\left( \cdot ,x,y \right) :\left( x,y \right) \in B_{\mathbb{H}\times \mathbb{H}}\left( R \right) \right\}, \cr
		\left\{ F_2\left( \cdot ,x,y \right) :\left( x,y \right) \in B_{\mathbb{H}\times \mathbb{H}}\left( R \right) \right\},\quad
		\left\{ G_2\left( \cdot ,x,y,z \right) :\left( x,y,z \right) \in B_{\mathbb{H}\times \mathbb{H}}\left( R \right) \times \mathbb{Z} \right\}, \nonumber
		\end{gather}
		are uniformly almost periodic.}
\end{rem}

As we know above, $ A_2\left( \cdot \right) $ is periodic, \remref{rem4.5} holds and the family of measures  $\varLambda$ is tight. By proceeding as \cite{prato1995periodic}, we can prove that the mapping 
$$t \mapsto \mu _{t}^{x}, \qquad t\in \mathbb{R},\ x\in\mathbb{H} $$ 
is almost periodic. 

\section{The averaged equation}\label{sec-5}
\begin{lem}\label{lem6.1} Under {\rm (A1)-(A5)}, for any compact set  $ \mathbb{K}\subset \mathbb{H}, $ 
	the family of functions 
	\begin{eqnarray}\label{en61}
	\left\{ t\in \mathbb{R}\mapsto \int_{\mathbb{H}}{B_1\left(t, x,y \right)}\mu _{t}^{x}\left( dy \right), \ x\in \mathbb{K}\right\} 
	\end{eqnarray}
	is uniformly almost periodic.
\end{lem}
\para{Proof:} As $ \mathbb{K} $ is a compact set in $ \mathbb{H}, $ so it is bounded and there exist some $ R>0 $ such that $ \mathbb{K}\subset B_\mathbb{H}(R). $ That is, for any $ x\in \mathbb{K}, $ we have $ \left\| x\right\| \leq R. $

Now, let us define
$$
\varPhi \left( t,x \right) =\int_{\mathbb{H}}{B_1\left( t,x,y \right)}\mu _{t}^{x}\left( dy \right),\quad \left( t,x \right) \in \mathbb{R}\times \mathbb{H}.
$$
Then, for any $ t,\tau \in \mathbb{R},\ n\in \mathbb{N} $ and $ x\in \mathbb{K}, $ according to the assumption (A3) and equation (\ref{en427}), we can get
\begin{eqnarray}
|\Phi \left( t+\tau ,x \right) -\Phi \left( t,x \right) |&\le& \Big| \int_{\lVert y \rVert_\mathbb{H} \le n}{B_1\left( t+\tau ,x,y \right)}\mu _{t+\tau}^{x}\left( dy \right) -\int_{\lVert y \rVert_\mathbb{H} \le n}{B_1\left( t+\tau ,x,y \right)}\mu _{t}^{x}\left( dy \right) \Big|\cr
&&+\Big| \int_{\lVert y \rVert_\mathbb{H} >n}{B_1\left( t+\tau ,x,y \right)}\mu _{t+\tau}^{x}\left( dy \right) \Big|+\Big| \int_{\lVert y \rVert_\mathbb{H} >n}{B_1\left( t+\tau ,x,y \right)}\mu _{t}^{x}\left( dy \right) \Big|\cr
&&+\Big| \int_{\mathbb{H}}{B_1\left( t+\tau ,x,y \right)}\mu _{t}^{x}\left( dy \right) -\int_{\mathbb{H}}{B_1\left( t,x,y \right)}\mu _{t}^{x}\left( dy \right) \Big|\cr
&\le& \underset{x\in \mathbb{K},\lVert y \rVert_\mathbb{H} \leq n}{\text{sup}}\lVert B_1\left( t+\tau ,x,y \right) \rVert _{\mathbb{H}}\Big| \int_{\mathbb{H}}{\left( \mu _{t+\tau}^{x}-\mu _{t}^{x} \right)}\left( dy \right) \Big|+\frac{c\left( 1+\lVert x \rVert \right)}{n}\cr
&&+\Big| \int_{\mathbb{H}}{\lVert B_1\left( t+\tau ,x,y \right) -B_1\left( t,x,y \right) \rVert _{\mathbb{H}}}\mu _{t}^{x}\left( dy \right) \Big|,\nonumber
\end{eqnarray}
fixed some $ \bar{n} $ such that $ {c\left( 1+\lVert x \rVert \right)}/{\bar{n}}\leq \epsilon/3. $ Then, for any  $ t \in \mathbb{R} $, due to the mapping $t \mapsto \mu _{t}^{x} $ and the families of functions $ B_1\left(\cdot, x, y\right)  $ are almost periodic, we can find some $ \tau \in  \mathbb{R} $ such that  
\begin{eqnarray}
|\Phi \left( t+\tau ,x \right) -\Phi \left( t,x \right) |< \epsilon,\nonumber
\end{eqnarray}
so, the function $ \varPhi \left( \cdot ,x \right) $  is almost periodic for any  $ x\in \mathbb{K} $. 

By preceding as in \lemref{lem4.2}, we can get that under {\rm (A1)-(A4)}, for any fixed  $ x_1,x_2\in \mathbb{K} $ and $ p\geq 1 $, there exists $  c_p>0 $ such that
\begin{eqnarray}\label{en432}
\underset{s<t}{\sup}\mathbb{E}\left\| v^{x_1}\left( t;s,0 \right) -v^{x_2}\left( t;s,0 \right) \right\| _{\mathbb{H}}^{p}\le c_p\left\| x_1-x_2 \right\| _{\mathbb{H}}^{p},
\end{eqnarray}
according to (\ref{en414}) and (\ref{en415}), let $ s \rightarrow 0, $ it yields
\begin{eqnarray}\label{en433}
\underset{t\in \mathbb{R}}{\sup}\mathbb{E}\left\| \eta ^{x_1}\left(t \right) -\eta ^{x_2}\left(t \right)  \right\| _{\mathbb{H}}^{p}\le c_p\left\| x_1-x_2 \right\| _{\mathbb{H}}^{p}.
\end{eqnarray}
Hence, thanks to (\rm A3) and (\ref{en433}), we can conclude that for any $ x_1,x_2\in \mathbb{K}, $  have
\begin{eqnarray}
\left\| \varPhi \left( t,x_1 \right) - \varPhi \left( t,x_2 \right) \right\| _{\mathbb{H}}&\le& \mathbb{E}\left\| B_1\left(t, x_1,\eta ^{x_1}\left( t \right) \right) -B_1\left(t, x_2,\eta ^{x_2}\left( t \right) \right) \right\| _{\mathbb{H}}\cr
&\le& c\big( \left\| x_1-x_2 \right\| _{\mathbb{H}}+\big( \mathbb{E} \left\| \eta ^{x_1}\left( t \right) -\eta ^{x_2}\left( t \right) \right\|^2_{\mathbb{H}}\big)^\frac{1}{2} \big).\cr
&\le& c\left\| x_1-x_2 \right\| _{\mathbb{H}}.\nonumber
\end{eqnarray} 
This means that the family of functions $ \left\{ \varPhi \left( t,\cdot \right) :\ t\in \mathbb{R} \right\} $  is equicontinuous about $ x. $ 

In view of the Theorem 2.10 in \cite{fink1974almost}, the above conclusions indicate that $ \left\{ \varPhi \left( \cdot ,x \right) : x\in \mathbb{K} \right\} $ is uniformly almost periodic. Hence, the proof is complete.  \qed

According to Theorem 3.4 in \cite{cerrai2017averaging}, we define 
\begin{eqnarray}\label{en17}
\bar{B}_1\left( x \right) :=\underset{T\rightarrow \infty}{\lim}\frac{1}{T}\int_0^T{\int_{\mathbb{H}}{B_1\left( t,x,y \right)}\mu _{t}^{x}\left( dy \right) dt, \quad x\in \mathbb{H}},
\end{eqnarray}
and thanks to {\rm (A3)} and (\ref{en427}), we have that
\begin{eqnarray}\label{en62}
\left\| \bar{B}_1\left( x \right) \right\| _{\mathbb{H}}\le c\left( 1+\left\| x \right\| _{\mathbb{H}} \right).
\end{eqnarray}
\begin{lem}\label{lem6.2}Under {\rm (A1)-(A5)}, for any  $ T>0,s,\sigma\in \mathbb{R} $ and  $ x,y\in \mathbb{H} $,
	\begin{align}\label{en63}
	\mathbb{E}\Big\| \frac{1}{T}\int_s^{s+T}{B_1\left(t, x,v^x\left( t;s,y \right) \right)}dt-\bar{B}_1\left( x \right) \Big\| _{\mathbb{H}} \le \frac{c}{T }\left( 1+\left\| x \right\| _{\mathbb{H}} +\left\| y \right\| _{\mathbb{H}}  \right) +\alpha \left( T,x \right),
	\end{align}
	for some mapping  $ \alpha : \left[ 0,\infty \right) \times \mathbb{H}\rightarrow \left[ 0,\infty \right)  $ such that
	\begin{eqnarray}\label{en31}
	\underset{T>0}{\sup}\alpha \left( T,x \right) \le c\left( 1+\left\| x \right\| _{\mathbb{H}}  \right), \quad x\in \mathbb{H},
	\end{eqnarray}
	and for any compact set $ \mathbb{K}\subset\mathbb{H}, $ have 
	\begin{eqnarray}\label{en32}
	\underset{T\rightarrow \infty}{\lim}\underset{x\in \mathbb{K}}{\sup}\ \alpha \left( T,x \right) =0.
	\end{eqnarray}	
\end{lem}
\para{Proof:} We denote 
$$
\psi ^xB_1\left( t,y \right) :=B_1\left( t,x,y \right) -\int_{\mathbb{H}}{B_1\left(t, x,w \right)}\mu _{t}^{x}\left( dw \right).
$$
So 
\begin{eqnarray}\label{en64}
&&\mathbb{E}\bigg( \frac{1}{T}\int_s^{s+T}{\Big[ B_1\left(t, x,v^x\left( t;s,y \right) \right) -\int_{\mathbb{H}}{B_1\left(t, x,w \right)}\mu _{t}^{x}\left( dw \right) \Big]}dt \bigg) ^2 \cr
&=&\frac{2}{T^2}\int_s^{s+T}{\int_r^{s+T}{\mathbb{E}\left[ \psi ^xB_1\left( r,v^x\left( r;s,y \right) \right) \psi ^xB_1\left( t,v^x\left( t;s,y \right) \right) \right]}}dtdr \cr
&=&\frac{2}{T^2}\int_s^{s+T}{\int_r^{s+T}{\mathbb{E}\left[ \psi ^xB_1\left( r,v^x\left( r;s,y \right) \right) P_{r,t}^{x}\psi ^xB_1\left( r,v^x\left( r;s,y \right) \right) \right]}}dtdr \cr
&\leq&\frac{2}{T^2}\int_s^{s+T}{\int_r^{s+T}{\left( \mathbb{E} \left|  \psi ^xB_1\left( r,v^x\left( r;s,y \right) \right) \right| ^{2} \right) ^{\frac{1}{2}}\left( \mathbb{E} \left|  P_{r,t}^{x} \psi ^x B_1\left( r,v^x\left( r;s,y \right) \right) \right|^{2} \right)^{\frac{1}{2}}}}dtdr. \cr
&&
\end{eqnarray}
Due to {\rm (A3)}, (\ref{en413}) and (\ref{en427}), we have
\begin{eqnarray}
\mathbb{E}\left| \psi ^xB_1\left( r,v^x\left( r;s,y \right) \right)\right| ^{2} 
&\leq& c\mathbb{E}\left\| B_1\left(r, x,v^x\left( r;s,y \right) \right) \right\| _{\mathbb{H}}^{2}+c\mathbb{E}\Big( \int_{\mathbb{H}}{\left\| B_1\left(r, x,w \right) \right\| _{\mathbb{H}}}\mu _{r}^{x}\left( dw \right) \Big)^2 \cr
&\leq& c\left( 1+\left\| x \right\| _{\mathbb{H}}^{2}+\mathbb{E}\left\| v^x\left( r;s,y \right) \right\| _{\mathbb{H}}^{2} \right) \cr
&\leq& c\left( 1+\left\| x \right\| _{\mathbb{H}}^{2}+e^{-2\delta \left( r-s \right)}\left\| y \right\| _{\mathbb{H}}^{2} \right),\nonumber
\end{eqnarray}
and according to (\ref{en415}), we get
\begin{eqnarray}
\mathbb{E}\left|  P_{r,t}^{x}\psi ^xB_1\left( r,v^x\left( r;s,y \right) \right) \right| ^{2}
&=&\mathbb{E}\Big\| P_{r,t}^{x}\Big[ B_1\left(r, x,v^x\left( r;s,y \right) \right) -\int_{\mathbb{H}}{B_1\left(r, x,w \right)}\mu _{r}^{x}\left( dw \right) \Big] \Big\| _{\mathbb{H}}^{2}\cr
&=&\mathbb{E}\left\| B_1\left(t, x,v^x\left( t;s,y \right) -B_1\left( t,x,\eta ^x\left( t \right) \right) \right) \right\| _{\mathbb{H}}^{2}\cr
&\leq&c \mathbb{E}\left\| v^x\left( t;s,y \right) -\eta ^x\left( t \right) \right\| _{\mathbb{H}}^{2}\cr
&\leq& c e^{-\delta_2 \left( t-s \right)}\left( 1+\left\| x \right\| _{\mathbb{H}}^{2}+\left\| y \right\| _{\mathbb{H}}^{2} \right).\nonumber
\end{eqnarray} 
Let $ \delta ={\delta _2}/{2}$, it follows
\begin{eqnarray}\label{en66}
&&\mathbb{E} \Big\| \frac{1}{T}\int_s^{s+T}{\Big[ B_1\left(t, x,v^x\left( t;s,y \right) \right) -\int_{\mathbb{H}}{B_1\left(t, x,w \right)}\mu _{t}^{x}\left( dw \right) \Big]}dt \Big\| _{\mathbb{H}} \cr
&\leq& \frac{c }{T }\left( 1+\left\| x \right\| _{\mathbb{H}} +\left\| y \right\| _{\mathbb{H}} \right) \Big( \int_s^{s+T}{\int_r^{s+T}{e^{-\delta \left( t-s \right)}}}dtdr\Big)^{\frac{1}{2}}\cr
&\leq& \frac{c }{T}\left( 1+\left\| x \right\| _{\mathbb{H}} +\left\| y \right\| _{\mathbb{H}}  \right).
\end{eqnarray}

Since the family of functions (\ref{en61}) is uniformly almost periodic, according to Theorem 3.4 in \cite{cerrai2017averaging}, we can get that the limit
$$
\underset{T\rightarrow  \infty}{\lim}\frac{1}{T}\int_s^{s+T}{\int_{\mathbb{H}}{B_1\left(t, x,w \right)}}\mu _{t}^{x}\left( dw \right) dt 
$$
converges to  $ \bar{B}_1\left( x \right) $ uniformly with respect to $ s\in \mathbb{R} $ and $ x\in \mathbb{K} $ ($ \mathbb{K} $ is a compact set in  $ \mathbb{H} $).  Therefore, if we define
$$\alpha \left( T,x \right) = \Big\| \frac{1}{T}\int_s^{s+T}{\int_{\mathbb{H}}{B_1\left(t,x,w \right)}\mu _{t}^{x}\left( dw \right) dt}-\bar{B}_1\left( x \right) \Big\| _{\mathbb{H}}, $$
we get the conclusion (\ref{en32}). Moreover, 
thanks to the assumption (A3), equation (\ref{en427}) and (\ref{en62}), we also can get
$$
\alpha \left( T,x \right) \le  \Big| \frac{1}{T}\int_s^{s+T}{\int_{\mathbb{H}}{\lVert B_1\left( t,x,w \right) \rVert _{\mathbb{H}}}\mu _{t}^{x}\left( dw \right) dt} \Big|+\lVert \bar{B}_1\left( x \right) \rVert _{\mathbb{H}}\le c\left( 1+\left\| x \right\| _{\mathbb{H}}  \right).
$$
Hence, the proof is complete.
\qed

Now, we introduce the averaged equation
\begin{eqnarray}\label{en67}
du\left( t \right) &=&\left[ A_1(t)u\left( t \right) +\bar{B}_1\left( u\left( t \right) \right) \right] dt+F_1\left(t, u\left( t \right) \right) dw^{Q_1}\left( t \right) \cr
&&+\int_{\mathbb{Z}}{G_1\left( t,u\left( t \right) ,z \right)}\tilde{N}_1\left( dt,dz \right), \quad\quad  u\left( 0 \right) =x\in \mathbb{H}.
\end{eqnarray}
Due to the assumption {\rm (A3)}, we can easily get that the mapping $ \bar{B}_1 $  is Lipschitz continuous. So, for any  $ x\in \mathbb{H},  T>0 $ and  $ p\ge 1 $, equation  (\ref{en67}) admits a unique mild solution $\bar{u}\in L^p\left( \varOmega ;D\left( \left[ 0,T \right] ;\mathbb{H} \right) \right). $

\section{Averaging principles}\label{sec-6}
In this section, we prove that the slow motion $ u_{\epsilon}  $ converges to the averaged motion  $ \bar{u} $, as $ \epsilon \rightarrow 0 $.
\begin{thm}\label{thm7.1}Under {\rm (A1)-(A5)}, fix $ x\in \mathcal{D}(( -A_1 ) ^{\theta}) \left( \theta \in [ 0,\bar{\theta} ) \right) $ and  $ y\in \mathbb{H} $, for any  $ T>0 $ and  $ \eta >0 $, we have 
	\begin{eqnarray}\label{en71}
	\underset{\epsilon \rightarrow 0}{\lim}\mathbb{P}\Big( \underset{t\in \left[ 0,T \right]}{\sup}\left\| u_{\epsilon}\left( t \right) -\bar{u}\left( t \right) \right\| _{\mathbb{H}}>\eta \Big) =0,
	\end{eqnarray}
	where  $ \bar{u} $ is the solution of the averaged equation (\ref{en67}).
\end{thm}
\para{Proof:} For any  $ h\in \mathcal{D}\left( A_1 \right) \cap L^{\infty}\left(\mathcal{O} \right) $, we have
\begin{eqnarray}
\left< u_{\epsilon}\left( t \right) ,h \right> _{\mathbb{H}}
&=&\left< x,h \right> _{\mathbb{H}}+\int_0^t{\left< A_1(r) u_{\epsilon}\left( r \right) ,h \right> _{\mathbb{H}}}dr+\int_0^t{\left< \bar{B}_1\left( u_{\epsilon}\left( r \right) \right) ,h \right> _{\mathbb{H}}}dr\cr
&&+{\Big< \int_0^t F_1\left(r, u_{\epsilon}\left( r \right) \right) dw^{Q_1}\left( r \right),h\Big> _{\mathbb{H}}}\cr
&&+{\Big<\int_0^t \int_{\mathbb{Z}}{G_1\left( r,u_\epsilon\left( r \right) ,z \right)}\tilde{N}_1\left( dr,dz \right) ,h \Big> _{\mathbb{H}}}+R_{\epsilon}\left( t \right),\nonumber
\end{eqnarray} 
where
$$
R_{\epsilon}\left( t \right) :=\int_0^t{\left< B_1\left(r, u_{\epsilon}\left( r \right) ,v_{\epsilon}\left( r \right) \right) -\bar{B}_1\left( u_\epsilon\left( r \right) \right) ,h \right> _{\mathbb{H}}dr}.
$$
According to the proof in Section \ref{sec-3}, we know that the family $\left\{  \mathcal{L}\left(  u_{\epsilon} \right) \right\} _{\epsilon \in \left( 0,1 \right]}$ is tight in $ \mathcal{P} \left( D\left(\left[0,T \right]; \mathbb{H}\right)\right)$. Hence, in order to prove \thmref{thm7.1}, it is sufficient to prove  $ \underset{\epsilon \rightarrow 0}{\lim}\mathbb{E}\underset{t\in \left[ 0,T \right]}{\sup}\left| R_{\epsilon}\left( t \right) \right|=0 $.

For any  $ \epsilon >0 $ and some deterministic   constant  $ \delta _{\epsilon}>0 $, we divide the interval $ \left[ 0,T \right]  $ in subintervals of the size  $ \delta _{\epsilon} $. In each time interval  $ \left[ k\delta _{\epsilon},\left( k+1 \right) \delta _{\epsilon} \right] ,k=0,1,\cdots ,\lfloor  {T}/{\delta _{\epsilon}} \rfloor  $, we define the following auxiliary fast motion $ \hat{v}_{\epsilon} $ 
\begin{eqnarray}\label{en72}
d\hat{v}_{\epsilon}\left( t \right) &=&\frac{1}{\epsilon}\left[ \left( A_2\left( t \right) -\alpha \right) \hat{v}_{\epsilon}\left( t \right) +B_2\left( t,u_{\epsilon}\left( k\delta _{\epsilon} \right), \hat{v}_{\epsilon}\left( t \right) \right) \right] dt\cr
&&+\frac{1}{\sqrt{\epsilon}}F_2\left( t,u_{\epsilon}\left( k\delta _{\epsilon} \right), \hat{v}_{\epsilon}\left( t \right) \right) d\omega ^{Q_2}\left( t \right) \cr
&&+\int_{\mathbb{Z}}{G_2}\left( t,u_{\epsilon}\left( k\delta _{\epsilon} \right), \hat{v}_{\epsilon}\left( t \right) \right) \tilde{N}_{2}^{\epsilon}\left( dt,dz \right).
\end{eqnarray}
According to the definition of  $ \hat{v}_{\epsilon} $, we know that an analogous estimate to (\ref{en34}) holds. So, for any  $ p\ge 1 $, we have 
\begin{eqnarray}\label{en73}
\int_0^T{\mathbb{E}\left\| \hat{v}_{\epsilon}\left( t \right) \right\| _{\mathbb{H}}^{p}}dt\le c_{p,T}\left( 1+\left\| x \right\| _{\mathbb{H}}^{p}+\left\| y \right\| _{\mathbb{H}}^{p} \right).
\end{eqnarray}
\begin{lem}\label{lem7.2}Under {\rm (A1)-(A5)}, fix  $ x\in \mathcal{D}(( -A_1 ) ^{\theta})  \left( \theta \in [ 0,\bar{\theta} ) \right)$ and  $ y\in \mathbb{H} $, there exists a constant  $ \kappa >0 $, such that if
	$$
	\delta _{\epsilon}=\epsilon \ln ^{\epsilon ^{-\kappa}},
	$$
	we have
	\begin{eqnarray}\label{en74}
	\underset{\epsilon \rightarrow 0}{\lim}\underset{t\in \left[ 0,T \right]}{\sup}\mathbb{E}\left\| \hat{v}_{\epsilon}\left( t \right) -v_{\epsilon}\left( t \right) \right\| _{\mathbb{H}}^{p}=0.
	\end{eqnarray}
\end{lem}
\para{Proof:} Let $ \epsilon >0 $ be fixed. For $ k=0,1,\cdots ,\lfloor {T}/{\delta _{\epsilon}} \rfloor $  and  $ t\in \left[ k\delta _{\epsilon},\left( k+1 \right) \delta _{\epsilon} \right]  $, let $ \rho_{\epsilon}\left( t \right)  $  be the solution of the problem
\begin{eqnarray}
d\rho _{\epsilon}\left( t \right) &=&\frac{1}{\epsilon}\left( A_2\left( t \right) -\alpha \right) \rho _{\epsilon}\left( t \right) dt+\frac{1}{\sqrt{\epsilon}}K_{\epsilon}\left( t \right) d\omega ^{Q_2}\left( t \right) \cr
&&+\int_{\mathbb{Z}}{H_{\epsilon}\left( t,z \right)}\tilde{N}_{2}^{\epsilon}\left( dt,dz \right), \quad\quad\qquad \rho_{\epsilon}\left( k\delta _{\epsilon} \right) =0,\nonumber
\end{eqnarray} 
where
$$
K_{\epsilon}\left( t \right) :=F_2\left( t,u_{\epsilon}\left( k\delta _{\epsilon} \right) ,\hat{v}_{\epsilon}\left( t \right) \right) -F_2\left( t,u_{\epsilon}\left(t \right) ,v_{\epsilon}\left( t \right) \right), 
$$
$$
H_{\epsilon}\left( t,z \right) :=G_2\left( t,u_{\epsilon}\left( k\delta _{\epsilon} \right) ,\hat{v}_{\epsilon}\left( t \right) ,z \right) -G_2\left( t,u_{\epsilon}\left( t \right) ,v_{\epsilon}\left( t \right) ,z \right). 
$$
We get 
\begin{eqnarray}\label{en1}
\rho _{\epsilon}\left( t \right) =\psi_{\alpha ,\epsilon,2}\left( \rho_{\epsilon};k\delta _{\epsilon} \right) \left( t \right) +\varGamma _{\epsilon}\left( t \right) +\varPsi _{\epsilon}\left( t \right), \quad t\in \left[ k\delta _{\epsilon},\left( k+1 \right) \delta _{\epsilon} \right],
\end{eqnarray}
where
$$
\varGamma _{\epsilon}\left( t \right) =\frac{1}{\sqrt{\epsilon}}\int_{k\delta _{\epsilon}}^t{U_{\alpha ,\epsilon,2}\left( t,r \right) K_{\epsilon}\left( r \right)}dw^{Q_2}\left( r \right), 
$$
$$
\varPsi _{\epsilon}\left( t \right) =\int_{k\delta _{\epsilon}}^t{\int_{\mathbb{Z}}{U_{\alpha ,\epsilon,2}\left( t,r \right) H_{\epsilon}\left( r,z \right)}}\tilde{N}_{2}^{\epsilon}\left( dr,dz \right). 
$$
If we denote $ \varLambda _{\epsilon}\left( t \right) :=\hat{v}_{\epsilon}\left( t \right) -v_{\epsilon}\left( t \right)  $  and $ \vartheta _{\epsilon}\left( t \right) := \varLambda _{\epsilon}\left( t \right) -  \rho _{\epsilon}\left( t \right)$, we have
$$
d\vartheta _{\epsilon}\left( t \right) =\frac{1}{\epsilon}\left[ \left( A_2\left( t \right) -\alpha \right) \vartheta _{\epsilon}\left( t \right) +B_2\left( t,u_{\epsilon}\left( k\delta _{\epsilon} \right) ,\hat{v}_{\epsilon}\left( t \right) \right) -B_2\left( t,u_{\epsilon}\left( t \right) ,v_{\epsilon}\left( t \right) \right) \right] dt.
$$
Due to the \lemref{lem3.3}, for $ \alpha>0 $ large enough, using Young's inequality, we have
\begin{eqnarray}
\frac{1}{p}\frac{d}{dt}\left\| \vartheta _{\epsilon}\left( t \right) \right\| _{\mathbb{H}}^{p}
&\leq& \frac{1}{\epsilon}\left< \left( \gamma_2\left( t \right) A_2+ L _2\left( t \right)-\alpha \right) \vartheta _{\epsilon}\left( t \right) ,\vartheta _{\epsilon}\left( t \right) \right> _{\mathbb{H}}\left\| \vartheta _{\epsilon}\left( t \right) \right\| _{\mathbb{H}}^{p-2}\cr
&&+\frac{1}{\epsilon}\left< B_2\left( t,u_{\epsilon}\left( k\delta _{\epsilon} \right) ,\hat{v}_{\epsilon}\left( t \right) \right) -B_2\left( t,u_{\epsilon}\left( t \right) ,v_{\epsilon}\left( t \right) \right) ,\vartheta _{\epsilon}\left( t \right) \right> _{\mathbb{H}}\left\| \vartheta _{\epsilon}\left( t \right) \right\| _{\mathbb{H}}^{p-2}\cr
&\leq& \frac{c}{\epsilon}\left\| \vartheta _{\epsilon}\left( t \right) \right\| _{\mathbb{H}}^{p}-\frac{\alpha}{\epsilon}\left\| \vartheta _{\epsilon}\left( t \right) \right\| _{\mathbb{H}}^{p}+\frac{c}{\epsilon}  \left\| u_{\epsilon}\left( k\delta _{\epsilon} \right) -u_{\epsilon}\left( t \right) \right\| _{\mathbb{H}}  \left\| \vartheta _{\epsilon}\left( t \right) \right\| _{\mathbb{H}}^{p-1}\cr
&&+\frac{c}{\epsilon} \left\| \hat{v}_{\epsilon}\left( t \right) -v_{\epsilon}\left( t \right) \right\| _{\mathbb{H}}   \left\| \vartheta _{\epsilon}\left( t \right) \right\| _{\mathbb{H}}^{p-1}\cr
&\leq& -\frac{\alpha}{2\epsilon}\left\| \vartheta _{\epsilon}\left( t \right) \right\| _{\mathbb{H}}^{p}+\frac{c_p}{\epsilon} \left\| u_{\epsilon}\left( k\delta _{\epsilon} \right) -u_{\epsilon}\left( t \right) \right\| _{\mathbb{H}}^{p}+\frac{c_p}{\epsilon}\left\| \hat{v}_{\epsilon}\left( t \right) -v_{\epsilon}\left( t \right) \right\| _{\mathbb{H}}^{p}\cr
&\leq& -\frac{\alpha}{2\epsilon}\left\| \vartheta _{\epsilon}\left( t \right) \right\| _{\mathbb{H}}^{p}+\frac{c_p}{\epsilon} \left( 1+\left\| x \right\| _{\theta}^{p}+\left\| y \right\| _{\mathbb{H}}^{p} \right) \big( \delta _{\epsilon}^{\beta  \left( \theta \right)p}+\delta _{\epsilon}\big) \cr
&&+\frac{c_p}{\epsilon}\left\| \hat{v}_{\epsilon}\left( t \right) -v_{\epsilon}\left( t \right) \right\| _{\mathbb{H}}^{p}.\nonumber
\end{eqnarray} 
Using the Gronwall inequality, we get
\begin{align}\label{en201}
\left\| \vartheta _{\epsilon}\left( t \right) \right\| _{\mathbb{H}}^{p}\leq \frac{c_p}{\epsilon} \left( 1+\left\| x \right\| _{\theta}^{p}+\left\| y \right\| _{\mathbb{H}}^{p} \right)  \big( \delta _{\epsilon}^{\beta \left( \theta \right)p+1}+\delta _{\epsilon}^2 \big)  +\frac{c_p}{\epsilon}\int_{k\delta _{\epsilon}}^{t}{\left\| \hat{v}_{\epsilon}\left( r \right) -v_{\epsilon}\left( r \right) \right\| _{\mathbb{H}}^{p}}dr.
\end{align}
By proceeding as Lemma 6.3 in \cite{cerrai2011averaging}, we prove that 
\begin{align}\label{en202}
\mathbb{E}\left\| \varGamma _{\epsilon}\left( t \right) \right\| _{\mathbb{H}}^{p} \leq \frac{c_p}{\epsilon}\left( 1+\left\| x \right\| _{\theta}^{p}+\left\| y \right\| _{\mathbb{H}}^{p} \right)  \big( \delta _{\epsilon}^{\beta \left( \theta \right)p+1}+\delta _{\epsilon}^2 \big)+\frac{c_p}{\epsilon}\int_{k\delta _{\epsilon}}^{t}{\mathbb{E}\left\| \hat{v}_{\epsilon}\left( r \right) -v_{\epsilon}\left( r \right) \right\| _{\mathbb{H}}^{p} }dr.
\end{align}
For the other stochastic term $ \varPsi _{\epsilon}\left( t \right) $, using  Kunita's first inequality and the H\"{o}lder inequality, we yield  
\begin{eqnarray}\label{en203}
\mathbb{E}\left\| \varPsi _{\epsilon}\left( t \right) \right\| _{\mathbb{H}}^{p}
&\leq& c_p\mathbb{E}\Big(\frac{1}{\epsilon}\int_{k\delta _{\epsilon}}^t{\int_{\mathbb{Z}}{\left\| e^{-\frac{\alpha}{\epsilon}\left( t-r \right)}e^{\frac{\gamma_2\left(t,r \right)  }{\epsilon}A_2}H_{\epsilon}\left( r,z \right)\right\| } _{\mathbb{H}}^{2}}v_2\left( dz \right) dr \Big)^{\frac{p}{2}}\cr
&&+\frac{c_{p}}{\epsilon}\mathbb{E}\int_{k\delta _{\epsilon}}^t{\int_{\mathbb{Z}}{\left\| e^{-\frac{\alpha}{\epsilon}\left( t-r \right)}e^{\frac{\gamma_2\left(t,r \right)  }{\epsilon}A_2}H_{\epsilon}\left( r,z \right)\right\| } _{\mathbb{H}}^{p}}v_2\left( dz \right) dr \cr
&\leq& \frac{c_{p}}{\epsilon^{\frac{p}{2}}}\mathbb{E}\Big(\int_{k\delta _{\epsilon}}^{t}{e^{-\frac{2\alpha}{\epsilon}\left( t-r\right) }  \big( \left\| u_{\epsilon}\left( k\delta _{\epsilon} \right) -u_{\epsilon}\left( r \right) \right\| _{\mathbb{H}}^2+\left\| \hat{v}_{\epsilon}\left( r \right) -v_{\epsilon}\left( r \right) \right\| _{\mathbb{H}}^{2} \big)dr}\Big)^{\frac{p}{2}}\cr
&&+\frac{c_{p}}{\epsilon}\int_{k\delta _{\epsilon}}^{t}{ \mathbb{E}\left\| u_{\epsilon}\left( k\delta _{\epsilon} \right) -u_{\epsilon}\left( r \right) \right\| _{\mathbb{H}}^p+\mathbb{E}\left\| \hat{v}_{\epsilon}\left( r \right) -v_{\epsilon}\left( r \right) \right\| _{\mathbb{H}}^{p} }dr\cr
&\le&\frac{c_p}{\epsilon ^{\frac{p}{2}}}\int_{k\delta _{\epsilon}}^t{\left( \mathbb{E}\lVert u_{\epsilon}\left( k\delta _{\epsilon} \right) -u_{\epsilon}\left( r \right) \rVert _{\mathbb{H}}^{p}+\mathbb{E}\lVert \hat{v}_{\epsilon}\left( r \right) -v_{\epsilon}\left( r \right) \rVert _{\mathbb{H}}^{p} \right)}dr\Big( \int_{k\delta _{\epsilon}}^t{e^{-\frac{2\alpha}{\epsilon}\frac{p}{p-2}\left( t-r \right)}dr} \Big) ^{\frac{p-2}{2}} \cr
&&+\frac{c_p}{\epsilon}\int_{k\delta _{\epsilon}}^t{\left( \mathbb{E}\lVert u_{\epsilon}\left( k\delta _{\epsilon} \right) -u_{\epsilon}\left( r \right) \rVert _{\mathbb{H}}^{p}+\mathbb{E}\lVert \hat{v}_{\epsilon}\left( r \right) -v_{\epsilon}\left( r \right) \rVert _{\mathbb{H}}^{p} \right)} dr\cr
&\le& \frac{c_p}{\epsilon}\int_{k\delta _{\epsilon}}^t{\left( \mathbb{E}\lVert u_{\epsilon}\left( k\delta _{\epsilon} \right) -u_{\epsilon}\left( r \right) \rVert _{\mathbb{H}}^{p}+\mathbb{E}\lVert \hat{v}_{\epsilon}\left( r \right) -v_{\epsilon}\left( r \right) \rVert _{\mathbb{H}}^{p} \right)}  dr\cr
&\le&  \frac{c_p}{\epsilon}\big( 1+\lVert x \rVert _{\theta}^{p}+\lVert y \rVert _{\mathbb{H}}^{p} \big) \big( \delta _{\epsilon}^{\beta \left( \theta \right) p+1}+\delta _{\epsilon}^{2} \big) +\frac{c_p}{\epsilon}\int_{k\delta _{\epsilon}}^t{\mathbb{E}\lVert \hat{v}_{\epsilon}\left( r \right) -v_{\epsilon}\left( r \right) \rVert _{\mathbb{H}}^{p}}dr.
\end{eqnarray}
According to the Lemma 2.4 in \cite{cerrai2017averaging} and equations (\ref{en1})-(\ref{en203}), we obtain
\begin{eqnarray}
\mathbb{E}\left\| \hat{v}_{\epsilon}\left( t \right) -v_{\epsilon}\left( t \right) \right\| _{\mathbb{H}}^{p} 
 \leq  \frac{c_p}{\epsilon} \left( 1+\left\| x \right\| _{\theta}^{p}+\left\| y \right\| _{\mathbb{H}}^{p} \right)  \big( \delta _{\epsilon}^{\beta \left( \theta \right)p+1}+\delta _{\epsilon}^2 \big)+ \frac{c_p}{\epsilon}  \int_{k\delta _{\epsilon}}^{t}{\mathbb{E} \left\| \hat{v}_{\epsilon}\left( r \right) -v_{\epsilon}\left( r \right) \right\| _{\mathbb{H}}^{p}}dr.\nonumber
\end{eqnarray}
From the Gronwall lemma, this means
\begin{eqnarray}
\mathbb{E}\left\| \hat{v}_{\epsilon}\left( t \right) -v_{\epsilon}\left( t \right) \right\| _{\mathbb{H}}^{p} \leq  \frac{c_p}{\epsilon} \big( \delta _{\epsilon}^{\beta \left( \theta \right)p+1}+\delta _{\epsilon}^2 \big)e^{\frac{c_p}{\epsilon} \delta _{\epsilon}} \left( 1+\left\| x \right\| _{\theta}^{p}+\left\| y \right\| _{\mathbb{H}}^{p} \right) .\nonumber
\end{eqnarray}
For  $ t\in \left[ 0,T \right]  $, selecting  $ \delta _{\epsilon}=\epsilon \ln ^{\epsilon ^{-\kappa}} $, then if we take  $\kappa <\frac{\beta \left( \theta \right) p}{\beta \left( \theta \right) p+1+c_p}\land \frac{1}{2+c_p}    $, we have (\ref{en74}).\qed
\begin{lem}\label{lem7.3}Under the same assumptions as in \thmref{thm7.1}, for any  $ T>0 $, we have
	\begin{eqnarray}\label{en75}
	\underset{\epsilon \rightarrow 0}{\lim}\mathbb{E}\underset{t\in \left[ 0,T \right]}{\sup}\left| R_{\epsilon}\left( t \right) \right|=0.
	\end{eqnarray}
\end{lem}
\para{Proof:} According to the definition of  $ \bar{B}_1 $, we get that the mapping $ \bar{B}_1:\mathbb{H}\rightarrow \mathbb{H} $ is Lipschitz continuous. Using assumption {\rm (A3)}, \lemref{lem3.3} and \lemref{lem7.2}, we have
\begin{eqnarray}\label{en105}
&&\underset{\epsilon \rightarrow 0}{\lim}\mathbb{E}\underset{t\in \left[ 0,T \right]}{\sup}\left| R_{\epsilon}\left( t \right) \right|\cr
&\leq& \underset{\epsilon \rightarrow 0}{\lim}\mathbb{E}\int_0^T{\big| \left< B_1\left( r,u_{\epsilon}\left( r \right) ,v_{\epsilon}\left( r \right) \right) -B_1\left( r,u_{\epsilon}\left( \lfloor  {r}/{\delta _{\epsilon}} \rfloor \delta _{\epsilon} \right) ,\hat{v}_{\epsilon}\left( r \right) \right) ,h \right> _{\mathbb{H}} \big|dr}\cr
&&+\underset{\epsilon \rightarrow 0}{\lim}\mathbb{E}\underset{t\in \left[ 0,T \right]}{\sup}\Big| \int_0^t{\left< B_1\left(r, u_{\epsilon}\left( \lfloor  {r}/{\delta _{\epsilon}} \rfloor \delta _{\epsilon} \right) ,\hat{v}_{\epsilon}\left( r \right) \right) -\bar{B}_1 \left( u_{\epsilon}\left( r \right) \right) ,h \right> _{\mathbb{H}}dr} \Big|\cr
&\leq& \underset{\epsilon \rightarrow 0}{\lim}c_T\left\| h \right\| _{\mathbb{H}}\Big[ \underset{r\in \left[ 0,T \right]}{\sup}\mathbb{E}\left\| u_{\epsilon}\left( r \right) -u_{\epsilon}\left( \lfloor  {r}/{\delta _{\epsilon}} \rfloor \delta _{\epsilon} \right) \right\| _{\mathbb{H}}+\underset{r\in \left[ 0,T \right]}{\sup}\mathbb{E}\left\| \hat{v}_{\epsilon}\left( r \right) -v_{\epsilon}\left( r \right) \right\| _{\mathbb{H}} \Big] \cr
&&+\underset{\epsilon \rightarrow 0}{\lim}\sum_{k=0}^{\lfloor  {T}/{\delta _{\epsilon}} \rfloor}{\mathbb{E}}\Big| \int_{k\delta _{\epsilon}}^{\left( k+1 \right) \delta _{\epsilon}}{\left< B_1\left(r, u_{\epsilon}\left( \lfloor  {r}/{\delta _{\epsilon}} \rfloor \delta _{\epsilon} \right) ,\hat{v}_{\epsilon}\left( r \right) \right) -\bar{B}_1\left( u_{\epsilon}\left( k\delta _{\epsilon} \right) \right) ,h \right> _{\mathbb{H}}dr} \Big|\cr
&&+\underset{\epsilon \rightarrow 0}{\lim}\left\| h \right\| _{\mathbb{H}}\sum_{k=0}^{\lfloor  {T}/{\delta _{\epsilon}} \rfloor}{\int_{k\delta _{\epsilon}}^{\left( k+1 \right) \delta _{\epsilon}}{\mathbb{E}\left\| \bar{B}_1\left( u_{\epsilon}\left( k\delta _{\epsilon} \right) \right) -\bar{B}_1\left( u_{\epsilon}\left( r \right) \right) \right\| _{\mathbb{H}}dr}}\cr
&\leq& \underset{\epsilon \rightarrow 0}{\lim}c_T\left\| h \right\| _{\mathbb{H}}\Big[\left( 1+\left\| x \right\| _{\theta}+\left\| y \right\| _{\mathbb{H}} \right)  \big(\delta _{\epsilon}^{\beta \left( \theta \right)}+\delta _{\epsilon} \big)+\underset{r\in \left[ 0,T \right]}{\sup}\mathbb{E}\left\| \hat{v}_{\epsilon}\left( r \right) -v_{\epsilon}\left( r \right) \right\| _{\mathbb{H}}\Big]\cr
&&+\underset{\epsilon \rightarrow 0}{\lim}\sum_{k=0}^{\lfloor  {T}/{\delta _{\epsilon}} \rfloor}{\mathbb{E}}\Big| \int_{k\delta _{\epsilon}}^{\left( k+1 \right) \delta _{\epsilon}}{\left< B_1\left( r,u_{\epsilon}\left(  k \delta _{\epsilon} \right) ,\hat{v}_{\epsilon}\left( r \right) \right) -\bar{B}_1\left( u_{\epsilon}\left( k\delta _{\epsilon} \right) \right) ,h \right> _{\mathbb{H}}dr} \Big|\cr
&&+\underset{\epsilon \rightarrow 0}{\lim}c_T\left\| h \right\| _{\mathbb{H}}\left( 1+\left\| x \right\| _{\theta}+\left\| y \right\| _{\mathbb{H}} \right) \big(\delta _{\epsilon}^{\beta \left( \theta \right)+1}+\delta _{\epsilon}^2 \big)\left( \lfloor  {T}/{\delta _{\epsilon}} \rfloor +1\right). \nonumber
\end{eqnarray} 
So, we have to show that
\begin{eqnarray}\label{en76}
\underset{\epsilon \rightarrow 0}{\lim}\sum_{k=0}^{\lfloor  {T}/{\delta _{\epsilon}} \rfloor}{\mathbb{E}\Big| \int_{k\delta _{\epsilon}}^{\left( k+1 \right) \delta _{\epsilon}}{\left< B_1\left(r, u_{\epsilon}\left( k \delta _{\epsilon} \right) ,\hat{v}_{\epsilon}\left( r \right) \right) -\bar{B}_1\left( u_{\epsilon}\left( k\delta _{\epsilon} \right) \right) ,h \right> _{\mathbb{H}}dr} \Big|}=0.
\end{eqnarray}
If we set  $ \zeta _{\epsilon}={\delta _{\epsilon}}/{\epsilon} $, we have
\begin{eqnarray}
&&\mathbb{E}\Big| \int_{k\delta _{\epsilon}}^{\left( k+1 \right) \delta _{\epsilon}}{\left< B_1\left(r, u_{\epsilon}\left(k \delta _{\epsilon} \right) ,\hat{v}_{\epsilon}\left( r \right) \right) -\bar{B}_1\left( u_{\epsilon}\left( k\delta _{\epsilon} \right) \right) ,h \right> _{\mathbb{H}}dr} \Big|\cr
&=&\mathbb{E}\Big| \int_0^{\delta _{\epsilon}}{\left< B_1\left(k\delta_\epsilon+r, u_{\epsilon}\left( k\delta _{\epsilon} \right) ,\hat{v}_{\epsilon}\left( k\delta _{\epsilon}+r \right) \right) -\bar{B}_1\left( u_{\epsilon}\left( k\delta _{\epsilon} \right) \right) ,h \right> _{\mathbb{H}}dr} \Big|\cr
&=&\mathbb{E}\Big| \int_0^{\delta _{\epsilon}}{\left< B_1\big(k\delta_\epsilon+r, u_{\epsilon}\left( k\delta _{\epsilon} \right) ,\tilde{v}^{u_{\epsilon}\left( k\delta _{\epsilon} \right) ,v_{\epsilon}\left( k\delta _{\epsilon} \right)}\left( {r}/{\epsilon} \right) \big) -\bar{B}_1\left( u_{\epsilon}\left( k\delta _{\epsilon} \right) \right) ,h \right> _{\mathbb{H}}dr} \Big|\cr
&=&\delta _{\epsilon}\mathbb{E}\Big| \frac{1}{\zeta _{\epsilon}}\int_0^{\zeta _{\epsilon}}{\left< B_1\big( k\delta_\epsilon+\epsilon r,u_{\epsilon}\left( k\delta _{\epsilon} \right) ,\tilde{v}^{u_{\epsilon}\left( k\delta _{\epsilon} \right) ,v_{\epsilon}\left( k\delta _{\epsilon} \right)}\left( r \right) \big) -\bar{B}_1\left( u_{\epsilon}\left( k\delta _{\epsilon} \right) \right) ,h \right> _{\mathbb{H}}dr} \Big|,\nonumber
\end{eqnarray} 
where $ \tilde{v}^{u_{\epsilon}\left( k\delta _{\epsilon} \right) ,v_{\epsilon}\left( k\delta _{\epsilon} \right)}\left( r \right) $  is the solution of the fast motion equation (\ref{en41}) with the initial
datum given by  $ v_{\epsilon}\left( k\delta _{\epsilon} \right) $ and the frozen slow component given by $ u_{\epsilon}\left( k\delta _{\epsilon} \right) $. In addition, the noises in (\ref{en41}) are independent of $ u_{\epsilon}\left( k\delta _{\epsilon} \right) $ and $ v_{\epsilon}\left( k\delta _{\epsilon} \right) $. According to the proof of \lemref{lem6.2}, we get
\begin{eqnarray}
&&\mathbb{E}\Big| \int_{k\delta _{\epsilon}}^{\left( k+1 \right) \delta _{\epsilon}}{\left< B_1\left(r, u_{\epsilon}\left( k \delta _{\epsilon} \right) ,\hat{v}_{\epsilon}\left( r \right) \right) -\bar{B}_1\left( u_{\epsilon}\left( k\delta _{\epsilon} \right) \right) ,h \right> _{\mathbb{H}}dr} \Big|\cr
&\leq&\delta _{\epsilon}\frac{c}{{\zeta _{\epsilon}}}\left( 1+\mathbb{E}\left\| u_{\epsilon}\left( k\delta _{\epsilon} \right) \right\| _{\mathbb{H}}+\mathbb{E}\left\| v_{\epsilon}\left( k\delta _{\epsilon} \right) \right\| _{\mathbb{H}} \right) \left\| h \right\| _{\mathbb{H}}+\delta _{\epsilon}\left\| h \right\| _{\mathbb{H}}  \mathbb{E}\alpha \left( \zeta _{\epsilon},u_{\epsilon}\left( k\delta _{\epsilon} \right) \right).\nonumber
\end{eqnarray} 
In Section \ref{sec-3}, we proved that the family
	$$
	\left\{ u_{\epsilon}\left( k\delta _{\epsilon} \right) : \epsilon >0, \   k=0,\cdots, \lfloor  {T}/{\delta _{\epsilon}} \rfloor \right\} 
	$$
	is tight. Then, by proceeding as the proof of the equation (8.21) in \cite{cerrai2017averaging}, we also can get that for any  $ \eta >0 $, there exists a compact set  $ \mathbb{K}_{\eta}\subset \mathbb{H} $, such that
	\begin{eqnarray}
	\mathbb{E}\alpha \left( \zeta _{\epsilon},u_{\epsilon}\left( k\delta _{\epsilon} \right) \right) 
	&\leq& \underset{x\in \mathbb{K}_{\eta}}{\sup}\alpha \left( \zeta _{\epsilon},x \right) +\sqrt{\eta}c\left( 1+\left\| x \right\| _{\mathbb{H}}+\left\| y \right\| _{\mathbb{H}} \right).\nonumber
	\end{eqnarray} 
	 Due to the arbitrariness of  $ \eta $ and according to the \lemref{lem6.2}, 
we can get equation (\ref{en76}). Furthermore, equation (\ref{en75}) holds. \qed

Through the above proof, \thmref{thm7.1} is established.\qed
\section*{Acknowledgments}
This work was partly supported by the National Natural Science Foundation of China under Grant No. 11772255, the Fundamental Research Funds for the Central Universities, the Research Funds for Interdisciplinary Subject of Northwestern Polytechnical University, the Shaanxi Project for Distinguished Young Scholars, the Shaanxi Provincial Key R\&D Program 2020KW-013 and 2019TD-010, and the Seed Foundation of Innovation and Creation for Graduate Students in Northwestern Polytechnical University (ZZ2018027).

\section*{References}
\bibliography{references}
\end{document}